\documentclass[11pt,a4paper]{article}
\usepackage[margin=28mm]{geometry}
\usepackage{fontspec}
\usepackage{amsmath,amssymb,amsthm,mathtools,booktabs,array,capt-of}
\usepackage{xcolor}
\usepackage[colorlinks=true,linkcolor=blue!50!black,citecolor=blue!50!black,urlcolor=blue!50!black]{hyperref}
\allowdisplaybreaks[2]
\newtheorem{theorem}{Theorem}[section]
\newtheorem{proposition}[theorem]{Proposition}
\newtheorem{lemma}[theorem]{Lemma}
\newcommand{\F}{\mathbb F}
\newcommand{\M}{\mathcal M}
\newcommand{\lam}{\lambda}
\newcommand{\ch}{\chi}
\newcommand{\eps}{\varepsilon}
\newcommand{\one}{\mathbf 1}
\DeclareMathOperator{\Tr}{Tr}
\DeclareMathOperator{\Nm}{N}
\DeclareMathOperator{\CT}{CT}
\title{Central coefficients of self-reciprocal irreducible polynomials over finite fields of odd order}
\hypersetup{pdftitle={Central coefficients of self-reciprocal irreducible polynomials over finite fields of odd order},
            pdfauthor={Hongfeng Wu and Li Zhu}}
\author{Hongfeng Wu$^{1}$ and Li Zhu$^{2}$
\setcounter{footnote}{-1}
\footnote{\texttt{whfmath@gmail.com} (H. Wu) and \texttt{izhumath@pku.edu.cn} (L. Zhu)}
\\[0.5ex]
\small $^{1}$College of Science, North China University of Technology, Beijing, China
\\
\small $^{2}$School of Mathematical Sciences, Guizhou Normal University, Guiyang, China}

\date{}

\begin{document}
\maketitle
\begin{abstract}
Let $q$ be an odd prime power. We determine which elements $a\in\F_q$
occur as the central coefficient of a monic self-reciprocal irreducible
polynomial of degree $2n$. Every $a$ occurs for $n\ge5$, whereas $1$
is missing in degree eight over $\F_3$. We determine all exceptions
in smaller degrees and express the nondegenerate quartic counts in terms
of elliptic-curve traces. If $N_q(n,a)$ denotes the number of such
polynomials, then $N_q(n,a)\sim q^{n-1}/(2n)$ as $q^n\to\infty$,
uniformly in $q$, $n\ge5$ and $a$. For fixed $q$ the error is
$O_q(q^{3n/4})$, and we give explicit lower bounds for $n\ge16$.
Writing $F(X)=X^n G(X+X^{-1})$ turns the central coefficient into a
linear form $\lam(G)$. For irreducible $G$, the condition that $F$ be
irreducible is that $G(2)G(-2)$ be a nonsquare. The pairing
$(A,B)\mapsto\lam(AB)$ has full rank on each space of polynomials
of bounded degree, giving character-sum estimates for the individual counts.
\end{abstract}

\section{Introduction}
For each $a\in\F_q$ and $n\ge1$, there are exactly $q^{n-1}$
monic self-reciprocal polynomials of degree $2n$ with central
coefficient $a$. Reciprocity pairs the coefficients on either side
of $X^n$, leaving $b_1,\ldots,b_{n-1}$ free in
\begin{equation}\label{eq:F}
 F(X)=X^{2n}+\sum_{i=1}^{n-1}b_i(X^{2n-i}+X^i)+aX^n+1,
 \qquad b_i\in\F_q.
\end{equation}
Let $N_q(n,a)$ count the irreducibles in this family. The quartic identity
\[
 X^4+b(X^3+X)+2X^2+1=(X^2+1)(X^2+bX+1)
\]
gives $N_q(2,2)=0$ over every finite field. In characteristic two
the paired terms cancel at $X=1$, leaving $F(1)=a$, so zero is
excluded in every even degree. In odd characteristic the quartic
example raises the question of how long a central coefficient can
force reducibility. Henceforth $q$ is odd, $\eta$ is the quadratic
character of $\F_q$ with $\eta(0)=0$, and $\eps=\eta(-1)$.
Integers used as coefficients denote their images in $\F_q$.

Wan~\cite{Wan} and Ham and Mullen~\cite{HamMullen} settled the
single-coefficient problem of Hansen and Mullen~\cite{HM} for general
monic irreducibles. For self-reciprocal irreducibles, counts with
prescribed leading coefficients~\cite{Gao22,Garefalakis11}
also apply to trailing coefficients by symmetry. A leading block
that reaches $X^n$ already determines all of $F$. For the coefficient
of $X^k$, Garefalakis and Kapetanakis~\cite{GK12,GK15}
classified the exceptions for $1\le k\le\lfloor n/2\rfloor$, and
Gao~\cite[Corollary~4(b),(c)]{Gao23} gave further existence criteria
for noncentral positions. At the central position $k=n$, the last
exception occurs in degree eight over $\F_3$.

\begin{theorem}[Existence]\label{thm:classification}
Let $q$ be an odd prime power. For $n\ge5$, every $a\in\F_q$ satisfies
$N_q(n,a)>0$. For $1\le n\le4$ the missing coefficients are as follows.
\begin{enumerate}
\item For $n=1$, $N_q(1,a)=0$ exactly when $\eta(a^2-4)\ne-1$.
\item For $n=2$, the missing set is $\{2,-2\}$ if $\eps=1$ and
      $\{2,6\}$ if $\eps=-1$.
\item For $n=3$, the only pair with $N_q(3,a)=0$ is $(q,a)=(3,0)$.
      For $n=4$, the only pair with $N_q(4,a)=0$ is $(q,a)=(3,1)$.
\end{enumerate}
\end{theorem}

For $n\ge5$ the question becomes how often each central coefficient
occurs. Summing $N_q(n,a)$ over $a$ recovers Carlitz's count by
degree~\cite{Carlitz}, which Yucas and Mullen~\cite{YucasMullen}
also obtained from the multiplicative orders of the roots. In particular,
\[
 \frac1q\sum_{a\in\F_q}N_q(n,a)
 =\frac{q^{n-1}}{2n}
  +O\!\left(\frac{q^{\lfloor n/2\rfloor-1}}{n}\right)
 \qquad(n\ge2).
\]
The proportion of irreducibles among all $q^n$ polynomials is therefore
asymptotic to $1/(2n)$. The next theorem gives the same proportion
among the $q^{n-1}$ polynomials with any one prescribed central coefficient.

\begin{theorem}[Distribution]\label{thm:distribution}
Uniformly over odd prime powers $q$, integers $n\ge5$ and $a\in\F_q$,
\[
 N_q(n,a)=\frac{q^{n-1}}{2n}\bigl(1+o(1)\bigr)
 \qquad(q^n\longrightarrow\infty).
\]
For fixed odd $q$, one has $N_q(n,a)=q^{n-1}/(2n)+O_q(q^{3n/4})$
as $n\to\infty$, uniformly in $a$. For every odd $q$ and $a\in\F_q$,
\[
 N_q(n,a)>\frac{q^{n-1}}{20n}\quad(n\ge16),\qquad
 N_q(n,a)>\frac{q^{n-1}}{10n}\quad(n\ge17).
\]
\end{theorem}

In Section~\ref{sec:reduction} we write $F(X)=X^nG(X+X^{-1})$
with $G$ monic of degree $n$. The central coefficient of $F$ is then
$\lam(G)$, the constant term of $G(X+X^{-1})$. For irreducible
$G$, Meyn's criterion~\cite{Meyn} states that $F$ is irreducible
precisely when $G(2)G(-2)$ is a nonsquare. Section~\ref{sec:low}
uses this criterion to count quartics and determine the exceptions in
degrees six and eight. Over $\F_3$ an irreducible $F$ must satisfy
$F(1)+(-1)^nF(-1)=0$, since $G(2)$ and $G(-2)$ are opposites.
Thus $a=b_1$ in degree six and $a=b_2+1$ in degree eight. The two
missing central values correspond to the outer-coefficient exceptions
of Garefalakis and Kapetanakis~\cite{GK15}.

On polynomials of degree less than $r$, the pairing $(A,B)\mapsto
\lam(AB)$ has rank $r$. Adding multiples of $A(2)B(2)$ and
$A(-2)B(-2)$ lowers the rank by at most two. Section~\ref{sec:analytic}
uses this bound to estimate sums over irreducibles. In
Section~\ref{sec:units} the nonvanishing of the two evaluations
sharpens the bilinear estimates. The cases $n=5,6$ in
Section~\ref{sec:fivesix} complete the distribution proof.
For existence, Section~\ref{sec:AS} gives an Artin--Schreier
construction when $n$ is the characteristic $p\ge5$, and
Section~\ref{sec:finite} supplies polynomials over the remaining small fields.

\section{Reduction to half the degree}\label{sec:reduction}
Let $\M_r$ be the set of monic polynomials of degree $r$ in $\F_q[Y]$,
with $\M_0=\{1\}$. Write $\F_q[Y]_{<r}$ for the polynomials of degree
less than $r$ and $\one_E$ for the indicator of $E$. The classical
substitution below~\cite[Section~2]{GK12} halves the degree. With
$\CT_X$ denoting the constant term of a Laurent polynomial, put
\begin{equation}\label{eq:lambda}
 \lam(G)=\CT_X G(X+X^{-1})
 =\sum_{j\ge0}\binom{2j}{j}[Y^{2j}]G,
 \qquad D=Y^2-4,\qquad \ch(G)=\eta(G(2)G(-2)).
\end{equation}
The function $\ch$ is multiplicative, and $\ch(G)=0$ exactly when
$\gcd(G,D)\ne1$. We call $n=\deg G$ the reduced degree.

\begin{lemma}\label{lem:transform}
For $n\ge1$, the map $G\mapsto F=X^nG(X+X^{-1})$ is a bijection from
$\M_n$ onto all monic self-reciprocal polynomials of degree $2n$.
For each $a\in\F_q$, its restriction to $\lam(G)=a$ gives precisely
the family in \eqref{eq:F}. The polynomial $F$ is irreducible if and
only if $G$ is irreducible and $\ch(G)=-1$.
\end{lemma}
\begin{proof}
The $n+1$ polynomials $X^n(X+X^{-1})^j$, $0\le j\le n$, have distinct
leading degrees $n+j$ and form a basis of the space of polynomials
satisfying $F(X)=X^{2n}F(1/X)$. Monicity fixes the coefficient of
the last basis element, and $[X^n]F=\lam(G)$. A factorization of $G$
induces one of $F$. If $G$ is irreducible with root $\beta$, then a
root $\alpha$ of $X^2-\beta X+1$ determines $\beta=\alpha+\alpha^{-1}$.
Its degree over $\F_q$ is $2n$ exactly when $\beta^2-4$ is a nonsquare
in $\F_{q^n}$. The norm of $\beta^2-4$ is $G(2)G(-2)$, and the quadratic character
on $\F_{q^n}$ is $\eta\circ\Nm_{\F_{q^n}/\F_q}$.
\end{proof}

For $n\ge2$, irreducible $G$ never vanishes at $\pm2$, and hence
\begin{equation}\label{eq:countG}
 N_q(n,a)=\frac12\sum_{\substack{G\in\M_n\\G\text{ irreducible}}}
 (1-\ch(G))\one_{\lam(G)=a}.
\end{equation}
Coefficientwise Frobenius and the substitution $X\mapsto-X$ give
\begin{equation}\label{eq:symmetry}
 N_q(n,a)=N_q(n,a^p)=N_q(n,(-1)^na),\qquad q=p^e.
\end{equation}

Let $\mu$ denote the integer or polynomial M\"obius function.
The polynomial von Mangoldt function satisfies $\Lambda(P^j)=\deg P$
for monic irreducible $P$ and $j\ge1$, and vanishes otherwise.
Write $\pi_q(n)$ for the number of monic irreducibles of degree $n$ and put
\[
 Q_n(q)=\sum_{j=1}^{\lfloor n/2\rfloor}q^j.
\]
\begin{lemma}\label{lem:zero}
For every $n\ge1$,
\begin{equation}\label{eq:zerofreq}
 \sum_{f\in\M_n}\Lambda(f)=q^n,\qquad
 \sum_{f\in\M_n}\Lambda(f)\ch(f)=-1,\qquad
 \sum_{f\in\M_n}\Lambda(f)\ch(f)^2=q^n-2.
\end{equation}
Consequently,
\begin{equation}\label{eq:primezero}
 |n\pi_q(n)-q^n|\le Q_n(q),\qquad
 \left|\sum_{\substack{P\in\M_n\\P\text{ irreducible}}}\ch(P)\right|
 \le\frac{1+Q_n(q)}n.
\end{equation}
\end{lemma}
\begin{proof}
The generating series for the constant weight $1$ and the weight
$\ch$ are $(1-qz)^{-1}$ and $1-z$. Indeed,
$\sum_b\eta(b^2-4)=-1$, and in degrees at least two the two
evaluations are independent. Logarithmic differentiation of the
Euler products gives the first two identities. For the third, omit
the two prime powers $(Y\pm2)^n$ from the first sum. Finally, the
total von Mangoldt weight of proper powers is at most
$\sum_{d\mid n,\,d<n}d\pi_q(d)\le Q_n(q)$, which gives
\eqref{eq:primezero}.
\end{proof}

For $n\ge2$, put
\begin{equation}\label{eq:doubled-count}
 W_n(a)=\sum_{f\in\M_n}\Lambda(f)(1-\ch(f))\one_{\lam(f)=a}.
\end{equation}
By \eqref{eq:countG} and Lemma~\ref{lem:zero},
$0\le W_n(a)-2nN_q(n,a)\le 2Q_n(q)$, and the mean of $W_n(a)$
over $a\in\F_q$ is $(q^n+1)/q$.

\section{The exceptional coefficients in low degrees}\label{sec:low}
In degree two the discriminant of $X^2+aX+1$ gives the criterion in
Theorem~\ref{thm:classification}. Over $\F_3$, the condition
$\ch(G)=-1$ forces $G(1)+G(-1)=0$. For $g_j=[Y^j]G$ and
$\deg G\le5$ this reads $g_0+g_2+g_4=0$, and \eqref{eq:lambda} gives
\begin{equation}\label{eq:ternary-obstruction}
 \lam(G)=g_0-g_2=g_4-g_0 .
\end{equation}
A monic cubic with $\lam(G)=0$ or a monic quartic with $\lam(G)=1$
therefore has the factor $Y$.

We use the quadratic character bounds of Weil \cite{Weil}. For
$P\in\F_q[t]$ of degree $d\ge1$ that is not a constant times a square,
\begin{equation}\label{eq:Weileven}
 \left|\sum_{t\in\F_q}\eta(P(t))\right|\le
 \begin{cases}(d-1)\sqrt q,&d\text{ odd},\\
               (d-2)\sqrt q+1,&d\text{ even}.
 \end{cases}
\end{equation}
For squarefree $P$, the smooth projective model of $z^2=P(t)$ has
genus $\lfloor(d-1)/2\rfloor$. It has one point at infinity for odd
$d$, and $1+\eta(\operatorname{lc}P)$ rational points there for even
$d$. If $P=cH^2R$ with $R$ squarefree, omitting the zeros of $H$
changes the sum for $cR$ by at most $\deg H$, while $\deg R=d-2\deg H$.

\subsection{Quartic reciprocal polynomials}
Let $P,Q\in\F_q[t]$ be squarefree quadratics with no common
root, with $P$ monic and $Q$ of leading coefficient $c$.
Suppose that at each rational zero of either polynomial the
other has nonzero square value. Then
\begin{equation}\label{eq:two-nonsquares}
 \#\{t\in\F_q:\eta(P(t))=\eta(Q(t))=-1\}=\frac{q+1-\tau}{4}>0,
\end{equation}
where $\tau$ is the Frobenius trace of the smooth projective model
of $w^2=P(t)Q(t)$, which has genus one. By the condition at rational
zeros, the indicator
is $(1-\eta(P))(1-\eta(Q))/4$. The two quadratic sums are $-1$ and $-\eta(c)$. Counting the
$1+\eta(c)$ rational points at infinity gives
$\sum_t\eta(P(t)Q(t))=-\tau-\eta(c)$. Summing the indicator gives
\eqref{eq:two-nonsquares}, which is positive by $|\tau|\le2\sqrt q$.

\begin{proposition}\label{prop:n2}
For $n=2$,
\[
 N_q(2,2)=0,\qquad
 N_q(2,-2)=\begin{cases}0&\eps=1,\\(q+1)/2&\eps=-1,\end{cases}
 \quad
 N_q(2,6)=\begin{cases}(q-1)/2&\eps=1,\\0&\eps=-1.\end{cases}
\]
For $a\notin\{2,-2,6\}$ set $u=4(a-2)$ and $v=(a+2)^2/4$.
Let $t_a$ be the Frobenius trace of the smooth projective model of
\[
 C_a:\quad w^2=(t^2-u)(t^2-v).
\]
Then
\begin{equation}\label{eq:n2count}
 N_q(2,a)=\frac{q+1-\eps t_a}{4}>0.
\end{equation}
\end{proposition}
\begin{proof}
For $G=Y^2+bY+a-2$, Lemma~\ref{lem:transform} gives the following
two discriminant tests, obtained by Carlitz~\cite[Theorem~1]{CarlitzQuartic}
over prime fields,
\[
 \eta(b^2-4(a-2))=-1,\qquad
 \eta((a+2)^2-4b^2)=-1.
\]
At $a=2$ the first discriminant is $b^2$. At $a=-2$ the second test
requires $\eps=-1$ and $b\ne0$. In this case $b^2+16$ has no zero
and its character sum is $-1$, so precisely $(q+1)/2$ choices of $b$
pass the first test, all of them nonzero. At $a=6$ the tests are compatible
exactly when $\eps=1$. The polynomial $b^2-16$ has two zeros and
character sum $-1$, giving $(q-1)/2$ choices of $b$ in this case.

For $a\notin\{2,-2,6\}$, use $P=b^2-u$ and $Q=-4(b^2-v)$.
The parameters $u,v$ are nonzero, with $v-u=(a-6)^2/4$.
At $b^2=u$ the second discriminant is $(a-6)^2$,
and at $b^2=v$ the first is $(a-6)^2/4$. The curve $w^2=PQ$ is the
twist of $C_a$ by $-4$, with trace $\eps t_a$, so
\eqref{eq:two-nonsquares} gives \eqref{eq:n2count}.
\end{proof}

The same trace is obtained from $y^2=x(x-u)(x-v)$, the quadratic
twist by $u$ of the Legendre curve with parameter $v/u$. Indeed,
grouping the quartic sum by $x=t^2$ gives
$-1+\sum_x\eta(x(x-u)(x-v))$.

\subsection{Sextic reciprocal polynomials}
\begin{proposition}\label{prop:n3}
$N_q(3,a)=0$ if and only if $(q,a)=(3,0)$.
\end{proposition}
\begin{proof}
Write
\[
 G=Y^3+bY^2+cY+a-2b,
 \qquad\ch(G)=\eta((a+2b)^2-4(c+4)^2).
\]
Let $\mathcal B_a$ be the set of these polynomials with $\ch(G)=-1$,
and put $M=|\mathcal B_a|$ and
$R=\sum_{G\in\mathcal B_a}\#\{r\in\F_q:G(r)=0\}$.
The two factors $a+2b\pm2(c+4)$ are independent affine coordinates
in $(b,c)$. A reducible cubic has a root in $\F_q$, so
\[
 M=(q-1)^2/2,\qquad N_q(3,a)\ge M-R .
\]

Suppose first that the characteristic is not $3$. Put
$L_r=\eps\eta((r^2-1)(r^2-4))$ and
$S_q=\sum_r\eta((r^2-1)(r^2-4))=-1-t$, where $|t|\le2\sqrt q$
by the genus-one bound. No $G\in\mathcal B_a$ vanishes at $\pm2$.
For $r\ne\pm1,\pm2$, the number of members vanishing at $r$ is
\[
 \frac{q-2+L_r}{2}
       +\frac{1-qL_r}{2}\one_{a=2r}.
\]
Indeed, for $r\ne0$ substitute $c=-r^2-br+(2b-a)/r$ in $G(2)G(-2)$.
The two resulting linear factors in $b$ have product of leading
coefficients $-4(r^2-4)(r^2-1)/r^2$, and their zeros coincide
exactly when $a=2r$. The quadratic character sum gives the formula.
At $r=0$ use $b=a/2$. At $r=\pm1$ direct substitution gives
$(q-1)/2$ members if $a\ne2r$ and none otherwise. Thus
\begin{equation}\label{eq:cubicR}
 2R=q^2-4q+6+\eps S_q+
 \begin{cases}
 0 & a=\pm4,\\
 -(q-1)&a=\pm2,\\
 1-qL_{a/2}&\text{otherwise}.
 \end{cases}
\end{equation}
It follows that $M-R\ge(q-7-2\sqrt q)/2>0$ for $q\ge17$.
For $q=13$, $\eps=1$ gives $(q-5-2\sqrt q)/2>0$. For
$q=5,7,11$ the values of $S_q$ are $1,-1,3$, respectively.
Formula~\eqref{eq:cubicR} gives $M-R>0$ except at $(q,a)=(7,0)$,
where $M=R$. The cubic $Y(Y^2-1)\in\mathcal B_0$ contributes three
roots to $R$, so at most $R-2$ members are reducible in this case.

In characteristic $3$ the roots $r=\pm1$ equal $\mp2$ and contribute
zero. For any other $r$ the preceding calculation applies with
$L_r=\eps$, since $(r^2-1)(r^2-4)=(r^2-1)^2$. Hence
\[
 R=\frac{(q-2)(q-2+\eps)+(1-q\eps)\one_{a\notin\{\pm1\}}}{2}.
\]
This is less than $M$ except at $(q,a)=(3,0)$, already excluded
by \eqref{eq:ternary-obstruction}.
\end{proof}

\subsection{Octic reciprocal polynomials}
The first construction is Cohen's $R$-transform~\cite{CohenR}.
\begin{lemma}\label{lem:lift}
If $P$ is monic, self-reciprocal and irreducible of degree $2m$,
and $q^m\equiv1\pmod4$, then
\[
 \mathcal T(P)(X)=(2X)^{2m}P\left(\frac{X^2+1}{2X}\right)
\]
is monic, self-reciprocal and irreducible of degree $4m$.
\end{lemma}
\begin{proof}
Put $Q=q^m$ and let $\alpha$ be a root of $P$. Inversion is the
unique nonidentity involution on the cyclic Frobenius orbit that
commutes with Frobenius, so $\alpha^Q=\alpha^{-1}$. For
$z=\alpha-\alpha^{-1}\ne0$, we have $z^Q=-z$. Thus $z^2$ is a
nonsquare in $\F_Q$, as is
$\Nm_{\F_{Q^2}/\F_Q}(\alpha^2-1)=-z^2$ since $Q\equiv1\pmod4$.
The quadratic $X^2-2\alpha X+1$ is therefore irreducible over
$\F_{Q^2}$. Its root determines $\alpha$ and has degree $4m$
over $\F_q$. Expansion gives monicity and reciprocity.
\end{proof}

\begin{lemma}\label{lem:dquartic}
For $d\ne0$, the polynomial
$G=Y^4+tY^3+cY^2+dtY+d^2$ is irreducible if and only if
\[
 \eta(t^2-4(c-2d))=\eta((c+2d)^2-4dt^2)=-1.
\]
If $c\notin\{2d,-2d,6d\}$, such a $t$ exists.
\end{lemma}
\begin{proof}
Write $G=Y^2K(Y+d/Y)$, where $K=Z^2+tZ+c-2d$.
Its discriminant gives the first condition. For a root $\beta$ of
$K$, a root of $Y^2-\beta Y+d$ has degree four precisely when
$\beta^2-4d$ is a nonsquare in $\F_{q^2}$. Its norm is
$(c+2d)^2-4dt^2$, giving the second condition.
Under the additional hypothesis the two discriminants are
squarefree and have no common root. At a zero of the first, the
second equals $(c-6d)^2$. At a rational zero of the second, $d$ is a
square and the first equals $(c-6d)^2/(4d)$. Existence follows
from \eqref{eq:two-nonsquares}.
\end{proof}

For an irreducible monic $H$ of degree $j>1$ and
$\tau(Y)=(rY+s)/(tY+u)$ over $\F_q$ with $ru-st\ne0$, write $H[\tau]$ for the
monic normalization of $(tY+u)^jH(\tau(Y))$. Its leading coefficient
before normalization is $t^jH(r/t)$ if $t\ne0$ and $r^j$ otherwise.
It is irreducible of degree $j$ because $\tau^{-1}(\alpha)$ generates
the same field as a root $\alpha$ of $H$.

\begin{lemma}\label{lem:cayley}
If $J$ is monic irreducible of degree $m\ge1$, then $J(Z^2)$ is
irreducible exactly when $\eta((-1)^mJ(0))=-1$. In this case
\[
 \frac{(X-1)^{2m}}{J(1)}
 J\left(\left(\frac{X+1}{X-1}\right)^2\right)
\]
is monic, self-reciprocal and irreducible of degree $2m$.
For $m=4$, writing $J=T^4+bT^3+cT^2+dT+e$, its central coefficient is
\begin{equation}\label{eq:cayley4}
 \frac{70(1+e)-10(b+d)+6c}{1+e+b+d+c}.
\end{equation}
\end{lemma}
\begin{proof}
For a root $\beta$ of $J$, the square class of $\beta$ in
$\F_{q^m}$ is that of its norm $(-1)^mJ(0)$ in $\F_q$. This proves
the criterion for $J(Z^2)$. Apply the fractional linear substitution
$Z=(X+1)/(X-1)$ to this polynomial, as in the conjugated power maps
of~\cite{BassaMenares}. Inverting $X$ negates $Z$, giving reciprocity.
The normalization is valid since $J(1)=0$ would force $J=T-1$.
For $m=4$, the five central coefficients of
$(X+1)^{2j}(X-1)^{8-2j}$ are $70,-10,6,-10,70$ in increasing order
of $j$, which gives \eqref{eq:cayley4}.
\end{proof}
In particular, $(T+h)^4+C(T+h)^2+E$ is irreducible exactly when
$\eta(E)=\eta(C^2-4E)=-1$, by applying the square-root criterion
to $Z^2+CZ+E$ and translating.

\begin{proposition}\label{prop:n4}
$N_q(4,a)=0$ if and only if $(q,a)=(3,1)$.
\end{proposition}
\begin{proof}
Apply Lemma~\ref{lem:lift} to
$P=X^4+b(X^3+X)+AX^2+1$, since $q^2\equiv1\pmod4$. The central
coefficient of the lift is $8A+22$. Proposition~\ref{prop:n2}
therefore leaves only $38$ and either $6$ ($\eps=1$) or $70$
($\eps=-1$).

\smallskip\noindent
\emph{(i) The coefficient $6$ when $\eps=1$.}
Outside characteristic $3$, put $d=-4/3$ and $c=-8/9$ in
Lemma~\ref{lem:dquartic}. These parameters avoid its three exceptional
values of $c$. For an irreducible member of the family,
\[
 \lam(G)=6,\qquad
 G(2)G(-2)=-4(d+4)^2\bigl(t^2-4(c-2d)\bigr).
\]
The second expression is a nonsquare since $d+4\ne0$ and $\eps=1$.

Now let $q=3^{2k}$. Set
\[
 e=h-1,\quad E=h^4+h-1,\quad
 J=(T^2+2hT)^2+e=(T+h)^4-2h^2(T+h)^2+E.
\]
The conditions $\eta(e)=\eta(E)=-1$ imply the hypotheses of
Lemma~\ref{lem:cayley}, and \eqref{eq:cayley4} is zero, as required.
The quartic $E$ is squarefree since $E(-1)\ne0$ and $E'=h^3+1$.
Also $E(1)=1$, so $eE$ is squarefree.
At a zero of one factor the other is a square. The number of
suitable $h$ is therefore at least
$(q-6\sqrt q-1)/4>0$ for $q\ge81$.
For $q=9$, take $h^2=-1$ and
$J=(T+h)^4+(T+h)^2-1-h$.
Here $e=E=-1-h$ and both $E$ and $C^2-4E=h-1$ have norm $2$
to $\F_3$. Substitution of $b=d=h$, $c=1$ in \eqref{eq:cayley4}
gives zero.

\smallskip\noindent
\emph{(ii) The coefficient $70$ when $\eps=-1$.}
Outside characteristic $3$, put $d=-16,c=-96$ in
Lemma~\ref{lem:dquartic}. Then $\lam(G)=70$. The two discriminants
are $t^2+256$ and $64(t^2+256)$, and $G(2)G(-2)=64(196-9t^2)$.
Thus it suffices that
\[
 \eta(t^2+256)=-1,\qquad\eta(196-9t^2)=-1.
\]
Outside characteristics $5,7$, the equations for the roots are
$t^2=-256$ and $t^2=196/9$. These two values are nonzero and
differ by $(50/3)^2$, so the quadratics are squarefree and
have no common root. At a zero of one, the other is a nonzero square. Existence follows from \eqref{eq:two-nonsquares}.
Characteristic $5$ is incompatible with $\eps=-1$.
In characteristic $7$, the second condition holds for every
$t\ne0$, while the first has $(q+1)/2$ solutions, none equal to zero.

For $q=3^{2k+1}>3$, choose $h\notin\{0,1,-1\}$ and put
\[
 C=\frac{h(h^2+1)}{h+1},\qquad K=h^4+Ch^2,
 \qquad J=(T+h)^4+C(T+h)^2+E.
\]
The condition that \eqref{eq:cayley4} equal $1=70$ is
$2(b+d)+c=0$, which holds by the definition of $C$.
Moreover,
\[
 K=-\frac{h^3(h^2-h-1)}{h+1}\ne0,
 \qquad C^2+K=\frac{h^2(h-1)^2}{(h+1)^2}\ne0,
 \qquad C\ne0.
\]
Here $h^2-h-1$ has nonsquare discriminant $-1$.
It remains to choose $E$ with
$\eta(E)=\eta(C^2-E)=\eta(E+K)=-1$.
The three zeros $0,C^2,-K$ are distinct. At the first two, an
unvanished factor is $C^2$, and at the third it is $C^2+K$.
These are nonzero squares. The single character sums vanish and
the pairwise sums are $1,-1,1$. Thus the number
of choices is $(q+1-S)/8$, where
$S=\sum_E\eta(E(C^2-E)(E+K))$ and $|S|\le2\sqrt q$.
This number is positive.

\smallskip\noindent
\emph{(iii) The coefficient $38$.}
Equation \eqref{eq:cayley4} equals $38$ precisely when
$3(b+d)+2c=2(1+e)$. Both families below satisfy this identity.
\begin{align*}
\mathrm A:\quad&e_A=4h^2+6h-1,\quad E_A=h^4+e_A,
 &J_A&=(T^2+2hT)^2+e_A;\\
\mathrm B:\quad&E_B=3h^3+5h^2+6h-1,\quad D_B=h^4-4E_B,
 &J_B&=(T+h)^4-h^2(T+h)^2+E_B.
\end{align*}
For A, it suffices that $\eps=1$ and $\eta(e_A)=\eta(E_A)=-1$.
For B, it suffices that $\eta(E_B)=\eta(D_B)=-1$.

Suppose first that $\eps=1$ and the characteristic is not $13$.
The quadratic $e_A$ has discriminant $52\ne0$, so its character
sum is $-1$. If $E_A=(h^2+uh+v)^2$, then $u=0$, $v=2$, and
both $6=0$ and $5=0$ in $\F_q$, which is impossible. The two polynomials are
coprime since $E_A-e_A=h^4$ and $e_A(0)=-1$. By \eqref{eq:Weileven},
\[
 \left|\sum_h\eta(E_A)\right|\le2\sqrt q+1,
 \qquad\left|\sum_h\eta(e_AE_A)\right|\le4\sqrt q+1.
\]
At a zero of $e_AE_A$ the other factor is $h^4$ or $-h^4$, a
square. Hence A supplies at least $(q-6\sqrt q-1)/4>0$ parameters
when $q\ge49$.

Use B when $\eps=-1$ outside characteristic $3$, and in
characteristic $13$. These cases exclude characteristics $3,5$.
The relation $D_B+4E_B=h^4$ and $E_B(0)=-1$ give
$\gcd(E_B,D_B)=1$. If $D_B=(h^2+uh+v)^2$, coefficient comparison
would give $u=-6,v=-28$ and force the characteristic to divide
$\gcd(360,780)=60$. Consequently,
\[
 \left|\sum\eta(E_B)\right|\le2\sqrt q,\quad
 \left|\sum\eta(D_B)\right|\le2\sqrt q+1,\quad
 \left|\sum\eta(E_BD_B)\right|\le6\sqrt q.
\]
At a zero of one factor the other is a square. Thus B supplies
at least $(q-10\sqrt q-1)/4>0$ parameters when $q\ge103$.

For the remaining prime fields, Table~\ref{tab:octic} gives $h$ with
both $U,V$ nonsquares, where $(U,V)=(e_A,E_A)$ or $(E_B,D_B)$ in
the indicated family. The entries certify $U=rs^2$, $V=rt^2$.
\begin{center}
\begin{minipage}{\linewidth}
\centering
\setlength{\abovecaptionskip}{0pt}
\captionof{table}{Prime-field parameters for central coefficient $38$ in degree eight.}\label{tab:octic}
\smallskip
\begin{tabular}{rrrrrr@{\qquad}rrrrrr}
\toprule
$p$&$h$&family&$r$&$s$&$t$&$p$&$h$&family&$r$&$s$&$t$\\\midrule
5&2&A&2&1&2&37&5&A&2&3&9\\
7&1&B&-1&1&3&41&3&A&3&2&20\\
11&10&B&-1&4&1&43&4&B&-1&7&8\\
13&6&B&2&3&1&47&1&B&-1&9&2\\
17&5&A&3&3&6&59&1&B&-1&20&13\\
19&3&B&-1&3&4&67&10&B&-1&27&22\\
23&3&B&-1&8&10&71&4&B&-1&29&1\\
29&2&A&2&12&6&79&4&B&-1&10&23\\
31&1&B&-1&7&12&83&1&B&-1&30&36\\\bottomrule
\end{tabular}
\end{minipage}
\end{center}
The choices $r=-1,2,3$ occur for $p\equiv3\pmod4$,
$p\equiv5\pmod8$, and $p\equiv5\pmod{12}$, respectively.
Quadratic reciprocity gives $\eta(r)=-1$ in each row.
For family A the only nonprime fields below $49$ are $\F_9$ and
$\F_{25}$ with $\eps=1$. Write them as $\F_p(z)$, $z^2-z+2=0$.
Taking $h=z$ and $h=z+1$, respectively, gives $(e_A,E_A)=(z,z+2)$
and $(3z+1,z+4)$. All four entries have norm $2$ to $\F_p$. For
family B every nonprime field in characteristic $13$ has size at
least $169$. In the other eligible characteristics, $\eps=-1$
requires odd extension degree, so a nonprime field has size at least $7^3$.

Finally, if $q=3^e$ with $e$ odd, use $J_0=T^4+T+2$ over
$\F_3$. It has no rational root. A product of two irreducible
quadratics with its constant coefficient would have to be
$(T^2+1)(T^2\pm T+2)$, whose cubic coefficient is nonzero.
Thus $J_0$ is irreducible. Lemma~\ref{lem:cayley} gives an octic
with central coefficient $2=38$, which remains irreducible over
$\F_{3^e}$ because $\gcd(8,e)=1$.

The missing value over $\F_3$ follows from \eqref{eq:ternary-obstruction}.
\end{proof}

\section{Coefficient pairings and prime sums}\label{sec:analytic}
Define the monic Dickson polynomials $D_j$ by
$D_j(X+X^{-1})=X^j+X^{-j}$ for $j\ge1$, and define
$U_0=1$, $U_1=Y$, $U_{j+1}=YU_j-U_{j-1}$.

\begin{lemma}\label{lem:pairings}
For $r\ge1$, on $\F_q[Y]_{<r}$ the pairings
$(A,B)\mapsto\lam(AB)$ and $(A,B)\mapsto\lam(DAB)$
are nondegenerate. Their Gram determinants in the monomial basis
are $2^{r-1}$ and $(-2)^r$, respectively.
\end{lemma}
\begin{proof}
In the monic basis $1,D_1,\ldots,D_{r-1}$ the first Gram matrix is
$\operatorname{diag}(1,2,\ldots,2)$, by taking constant terms.
The identity
\[
 (X-X^{-1})U_j(X+X^{-1})=X^{j+1}-X^{-j-1}
\]
gives $\lam(DU_iU_j)=-2\delta_{ij}$. Both changes of basis have
determinant $1$.
\end{proof}
The even moments $\lam(Y^{2j})=\binom{2j}{j}$ yield the central
binomial Hankel determinants studied in~\cite{CiglerTyson}. For the
second pairing, $\lam(DY^{2j})=-2C_j$, where
$C_j=\binom{2j}{j}-\binom{2j}{j+1}$ is the $j$th Catalan number.

\subsection{Exponential sums over irreducible polynomials}
Fix the canonical additive character $\psi$ of $\F_q$. On
$\F_q((Y^{-1}))$ put $e(\theta)=\psi([Y^{-1}]\theta)$ and
$|Y|=q$. We use Hayes's estimate~\cite{Hayes} in the form given
by Pollack~\cite{Pollack}. For linear phases over irreducibles,
see also~\cite[Section~6.2]{Porritt}.

\begin{lemma}\label{lem:linearphase}
For $n\ge5$, $t\in\F_q^\times$ and $u,v\in\F_q$,
\begin{equation}\label{eq:linearphase}
 \left|\sum_{\substack{P\in\M_n\\P\text{ irreducible}}}
 \psi\bigl(t\lam(P)+uP(2)+vP(-2)\bigr)\right|
 \le K_n(q):=q^{n-\frac12\lfloor n/2\rfloor}
             +q^{n-\lceil n/2\rceil+2}.
\end{equation}
\end{lemma}
\begin{proof}
Put $\ell_j=t\lam(Y^j)+u2^j+v(-2)^j$ and $r=\lceil n/2\rceil$.
By Lemma~\ref{lem:pairings}, the Hankel matrix of $t\lam(Y^j)$
on indices $0,\ldots,r-1$ has rank $r$. Each evaluation adds a matrix
of rank at most one, so
\[
 \operatorname{rank}(\ell_{i+j})_{0\le i,j<r}\ge r-2\ge1 .
\]
For $\theta=\sum_{j=0}^{n-1}\ell_jY^{-j-1}$ the summand is
$\psi(\ell_n)e(\theta P)$. The rank bound ensures $\theta\ne0$.
Choose a reduced fraction $A/B$ with $B$ monic, $d=\deg B\le
\lfloor n/2\rfloor$ and $|\theta-A/B|<q^{-d-n/2}$ as
in~\cite[Lemma~3]{Pollack}. If $|\theta-A/B|\ge q^{-n}$,
\cite[Lemma~5]{Pollack} gives the first term of \eqref{eq:linearphase}.
Otherwise the first $n$ Laurent coefficients of $A/B$ agree with
those of $\theta$. Writing $B=\sum_{i=0}^db_iY^i$, $b_d=1$, gives
$\sum_{i=0}^db_i\ell_{j+i}=0$ for $0\le j<n-d$.
Since $2r-2\le n-1$, the Hankel matrix has rank at most $d$.
Thus $d\ge r-2$, and the additional term $q^{n-d}$ in
\cite[Lemma~6]{Pollack} is at most $q^{n-r+2}$.
\end{proof}

To impose the quadratic character in \eqref{eq:countG}, let
$\tau=\sum_x\eta(x)\psi(x)$. The Gauss sum identity and
$|\tau|=\sqrt q$ give
\[
 \ch(P)=\tau^{-2}\sum_{u,v\in\F_q^\times}\eta(uv)\psi(uP(2)+vP(-2)).
\]
For each $t\ne0$, Lemma~\ref{lem:linearphase} bounds
$\sum_P\psi(t\lam(P))$ by $K_n(q)$ and
$\sum_P\ch(P)\psi(t\lam(P))$ by $(q-1)^2K_n(q)/q$.
Apply additive orthogonality to \eqref{eq:countG} and use
\eqref{eq:primezero} at the zero frequency. This gives, for $n\ge5$,
\begin{equation}\label{eq:distributionexplicit}
 \left|N_q(n,a)-\frac{q^{n-1}}{2n}\right|
 \le \frac{q-1}{2q}\left(1+\frac{(q-1)^2}{q}\right)K_n(q)
      +\frac{1+2Q_n(q)}{2qn}.
\end{equation}
To compare the error with the main term, write $m=\lfloor n/2\rfloor$
and put
\[
 D_q(n)=\max_{a\in\F_q}\left|\frac{2nN_q(n,a)}{q^{n-1}}-1\right|,
 \qquad C_q=(q-1)\left(1+\frac{(q-1)^2}{q}\right).
\]
Dividing \eqref{eq:distributionexplicit} by $q^{n-1}/(2n)$ and
using $2m\le n\le2m+1$ and $Q_n(q)<q^{m+1}/(q-1)$ gives, for $n\ge5$,
\begin{equation}\label{eq:relative}
 D_q(n)\le(2m+1)C_q\bigl(q^{-m/2}+q^{2-m}\bigr)
          +q^{-2m}+\frac{2q^{1-m}}{q-1}.
\end{equation}
For $m\ge7$ the right side decreases with $m$, since
$(2m+3)/(2m+1)<\sqrt q$. It also decreases with $q\ge3$, since
\[
 \frac{qC'_q}{C_q}
 =\frac{q}{q-1}+\frac{q-q^{-1}}{q-1+q^{-1}}
 <\frac{2q}{q-1}\le3
\]
and both $m/2$ and $m-2$ exceed $3$. Substitution of
$(q,m)=(5,7),(5,8),(3,9)$, using $\sqrt5>11/5$ and $\sqrt3>5/3$,
therefore gives
\begin{equation}\label{eq:large-ranges}
 D_q(n)<
 \begin{cases}
  1    & q\ge5,\quad n\in\{14,15\},\\
  1/2  & q\ge5,\quad n\ge16,\\
  7/10 & q=3,\quad n\ge18 .
 \end{cases}
\end{equation}

\subsection{Bilinear estimates on coefficient fibres}

\begin{lemma}\label{lem:bilinear}
Let $r,s\ge0$, let $\psi$ be a nontrivial additive character of
$\F_q$, and let $A,B$ be complex sequences on $\M_r,\M_s$,
respectively. With $\|A\|_2^2=\sum_{f\in\M_r}|A(f)|^2$ and
likewise for $B$, we have
\[
 \left|\sum_{f\in\M_r,\,g\in\M_s}A(f)B(g)\psi(\lam(fg))\right|
 \le q^{\max(r,s)/2}\|A\|_2\|B\|_2.
\]
\end{lemma}
\begin{proof}
Write $f=Y^r+f_0$ and $g=Y^s+g_0$. The terms involving a leading
monomial multiply rows or columns by scalars of modulus $1$. If
$r\le s$, orthogonality and Lemma~\ref{lem:pairings} show that the
matrix $(\psi(\lam(f_0g_0)))$ times its conjugate transpose is
$q^sI$. Its operator norm is $q^{s/2}$. The other case is symmetric.
\end{proof}

For $h\in\M_m$, $m\ge2$, write $h=Dk+h_0$ with $k\in\M_{m-2}$,
$\deg h_0<2$. Interpolating $h_0$ from $u=h(2)$ and $v=h(-2)$ gives
\begin{equation}\label{eq:evaluation-line}
 \lam(bh)=\lam(Dbk)+\alpha_bu+\beta_bv,\qquad
 \alpha_b=\frac{\lam(b)}2+\frac{\lam(Yb)}4,\quad
 \beta_b=\frac{\lam(b)}2-\frac{\lam(Yb)}4 .
\end{equation}
For a nontrivial additive character $\psi$, write
\[
 S_\theta(d;m)=\sum_{h\in\M_m}\theta(dh)\psi(\lam(dh)).
\]
\begin{lemma}\label{lem:mass}
For $m\ge1$,
\[
 \sum_{d\in\M_r}|S_1(d;m)|\le
 \begin{cases}0&m\ge r+1,\\q^r&m\le r.\end{cases}
\]
For $m\ge2$,
\[
 \sum_{d\in\M_r}|S_\ch(d;m)|\le
 \begin{cases}0&m\ge r+3,\\q^{r+1}&m\le r+2.\end{cases}
\]
\end{lemma}
\begin{proof}
For fixed $d$, the untwisted inner sum vanishes unless $\lam(dv)=0$
for every $v\in\F_q[Y]_{<m}$.
When $m\le r$, these are $m$ independent conditions on the $r$
variable coefficients of $d$, leaving $q^{r-m}$ choices. Each
inner sum has modulus at most $q^m$. When $m>r$, nondegeneracy
excludes a monic $d$.

For the twist $\ch$, summing over $k$ in \eqref{eq:evaluation-line}
forces $\lam(Ddw)=0$ for all $\deg w<m-2$. By the second pairing
there are $q^{r-m+2}$ possible $d$ if $m\le r+2$, and none otherwise.
The sums over $u,v$ are quadratic Gauss sums. Their product vanishes
if $\alpha_d\beta_d=0$ and otherwise has modulus $q$. Thus each $d$
contributes at most $q^{m-1}$.
\end{proof}

We use Vaughan's identity for monic polynomials~\cite{Merai}. For
nonnegative integers $U,V$ and $\deg f>V$ it reads
\begin{align}\label{eq:Vaughan}
 \Lambda(f)={}&\sum_{\substack{d\mid f\\\deg d\le U}}\mu(d)\deg(f/d)
 -\sum_{\substack{de\mid f\\\deg d\le U,\,\deg e\le V}}
       \mu(d)\Lambda(e)
 +\sum_{\substack{de\mid f\\\deg d>U,\,\deg e>V}}
       \mu(d)\Lambda(e).
\end{align}
Indeed, if $M=\mu_{\le U}$ and $L=\Lambda_{\le V}$ denote truncations
by degree, then $\mu*1=\one_{f=1}$ gives
$\Lambda=M*(1*\Lambda)-M*L*1+(\mu-M)*(\Lambda-L)*1+L$,
and $(1*\Lambda)(f)=\deg f$.
Write $B_V(g)=\sum_{e\mid g,\,\deg e>V}\Lambda(e)$.

For $n\ge7$ and an integer $2\le U\le(n-3)/2$, put
\begin{align*}
 B(n,q,U)&=\sum_{r=U+1}^{n-U-1}
 \sqrt{(n-r)(n-r-U)(1-q^{-1})}\,q^{n-\min(r,n-r)/2},\\
 E_0(n,q,U)&=\sum_{r=\lceil n/2\rceil}^{2U}rq^r,\\
 E_1(n,q,U)&=\sum_{r=\lceil(n-2)/2\rceil}^{2U}rq^{r+1}.
\end{align*}
An empty sum is zero.
\begin{proposition}\label{thm:explicit}
With this notation,
\begin{equation}\label{eq:explicit}
 \left|N_q(n,a)-\frac{q^n+1}{2qn}\right|
 \le\frac{q-1}{2qn}(2B+E_0+E_1)+\frac{Q_n(q)}n.
\end{equation}
\end{proposition}
\begin{proof}
Apply \eqref{eq:Vaughan} with $V=U$ to
$T_\theta=\sum_{f\in\M_n}\Lambda(f)\theta(f)\psi(\lam(f))$,
where $\theta=1$ or $\ch$. Its first term vanishes by
Lemma~\ref{lem:mass}, since $n\ge2U+3$.
In the second term, combine $de=b$. Its coefficient has absolute
value at most $\sum_{e\mid b}\Lambda(e)=\deg b$. The same lemma
therefore bounds this term by $E_0$ or $E_1$.

In the third term $s=\deg g>U$. Each prime power $e$ of degree
$j>U$ divides $q^{s-j}$ members of $\M_s$, so \eqref{eq:zerofreq}
gives $\sum_{g\in\M_s}B_U(g)=(s-U)q^s$. Since $0\le B_U(g)\le s$,
its squared norm is at most $s(s-U)q^s$. The squared norm of
$\mu$ on $\M_r$ equals $(1-q^{-1})q^r$ for $r\ge2$, by the Euler
product for squarefree polynomials. Lemma~\ref{lem:bilinear}
therefore gives $|T_1|\le B+E_0$ and $|T_\ch|\le B+E_1$.

Additive orthogonality bounds the deviation of $W_n(a)$ from
its mean by $(q-1)(2B+E_0+E_1)/q$. The comparison with $2nN_q(n,a)$
after \eqref{eq:doubled-count} gives \eqref{eq:explicit}.
\end{proof}

For $7\le n\le9$, taking $U=2$ in Proposition~\ref{thm:explicit} gives
\begin{equation}\label{eq:relative-seven-nine}
 D_q(n)\ll q^{-1/2}.
\end{equation}
Indeed, $\min(r,n-r)\ge3$ gives $B\ll q^{n-3/2}$, and $r\le4$
gives $E_0+E_1\ll q^5$. The change of main term and the prime powers
add at most $q^{-n}+2q^{1-n}Q_n(q)$ to the relative error.

\section{Bilinear estimates with nonzero evaluations}\label{sec:units}
Set $\rho=\ch^2$, the indicator of coprimality to $D$. For every
irreducible of degree at least two, $\rho=1$, so the weight
$1-\ch$ in \eqref{eq:countG} can be replaced by $\rho-\ch$.
This excludes factors vanishing at $2$ or $-2$ from both sums.

\subsection{Character sums and moments}
\begin{lemma}\label{lem:unitmass}
For $m\ge2$, $\theta\in\{\ch,\rho\}$, and $L=r-m+2$, the sum
$\sum_{d\in\M_r}|S_\theta(d;m)|$ is at most $H_\theta(r,m)$,
where
\[
\begin{array}{c|cc}
L&H_\ch(r,m)&H_\rho(r,m)\\\hline
L<0&0&0\\
0&q^{r+1}&q^r\\
1&q^r(q-2)&q^{r-1}(3q-4)\\
2,3&q^{r+1}(1-q^{-1})^2&4q^r(1-q^{-1})^2\\
L\ge4&q^{r+1}(1-q^{-1})^4&4q^r(1-q^{-1})^4.
\end{array}
\]
\end{lemma}
\begin{proof}
Write $h=Dk+b$, $\deg b<2$, as in \eqref{eq:evaluation-line}.
The $D$-pairing in
Lemma~\ref{lem:pairings} restricts the support to
\[
 d=\sum_{j=m-2}^r c_jU_j,\qquad c_r=1.
\]
For $L\ge0$ define $C(T)=\sum c_jT^j=T^{m-2}A(T)$ with
$A\in\M_L$. The constant terms in $d$ and $Yd$ and the
identities $U_j(\pm2)=(\pm1)^j(j+1)$ give
\begin{align*}
 \lam(db)&=\tfrac12C(1)b(2)+\tfrac12C(-1)b(-2),\\
 d(2)&=C(1)+C'(1),&d(-2)&=C(-1)-C'(-1).
\end{align*}
Write $x=C(1)$, $y=C(-1)$, $u=d(2)$ and $v=d(-2)$. Both twists
vanish unless $uv\ne0$. In that case the two Gauss sums for $\ch$
give modulus $q^{m-1}\one_{xy\ne0}$. For $\rho$, the corresponding
modulus is $q^{m-2} I(x)I(y)$, where $I(0)=q-1$
and $I(z)=1$ for $z\ne0$. Thus the exact sum of the two masses is
\begin{equation}\label{eq:combined-mass}
 q^{m-2}\sum_{C:\,uv\ne0}\bigl(q\one_{xy\ne0}+I(x)I(y)\bigr).
\end{equation}
For $L\ge4$, multiplication by $T^{m-2}$ is invertible modulo
$(T^2-1)^2$. Hermite interpolation for the free coefficients of $A$
therefore makes $(x,y,u,v)$ run through $\F_q^4$ with multiplicity
$q^{L-4}$. Summing and using $\sum_z I(z)=2(q-1)$ gives the last
row. For $L=2,3$, $(x,y)$ are still independent and we drop the
condition on $uv$. For $L=1$, $A=T+c$, so exactly two choices
$c=\pm1$ annihilate one of $x,y$. For $L=0$ only $d=U_r$ can contribute, and both
sums vanish if $r+1=0$ in $\F_q$, since then $d(2)=0$.
For $L<0$ the support is empty. The same formulas hold for each
nontrivial character $z\mapsto\psi(tz)$.
\end{proof}

For monic polynomials coprime to $D$, put
\[
 V_j=\sum_{f\in\M_j}\rho(f)=
 \begin{cases}1&j=0,\\q-2&j=1,\\(q-1)^2q^{j-2}&j\ge2,\end{cases}
 \qquad A_r=\sum_{d\in\M_r}\mu(d)^2\rho(d).
\]
The Euler product gives
\begin{equation}\label{eq:Ar}
 \sum_{r\ge0}A_rz^r=\frac{1-qz^2}{(1-qz)(1+z)^2}.
\end{equation}
The number of degree-$j$ monic irreducibles coprime to $D$ is
$I_q^*(j)=j^{-1}\sum_{d\mid j}\mu(d)q^{j/d}-2\one_{j=1}$.

\begin{lemma}\label{lem:moment}
For $s>V$, write $l_j=\lfloor V/j\rfloor$. The first moment over
polynomials coprime to $D$ is
\[
 S_s(V)=\sum_{g\in\M_s}B_V(g)\rho(g)
 =\sum_{j=1}^s jI_q^*(j)\sum_{k=l_j+1}^{\lfloor s/j\rfloor}V_{s-kj}.
\]
If $s\le2V+1$, the second moment is
\begin{equation}\label{eq:moment}
 M_s(V)=\sum_{g\in\M_s}B_V(g)^2\rho(g)
 =\sum_{j=1}^s j^2I_q^*(j)
 \sum_{k=l_j+1}^{\lfloor s/j\rfloor}\bigl(2(k-l_j)-1\bigr)V_{s-kj}.
\end{equation}
\end{lemma}
\begin{proof}
A prime $P$ of degree $j$ and valuation $t$ in $g$ contributes
$j\max\{t-l_j,0\}$ to $B_V(g)$. For $P\nmid D$, the divisibility
condition $P^k\mid g$ holds for $V_{s-kj}$ polynomials of degree $s$
coprime to $D$. Summing over $k>l_j$ gives the first moment.
If $s\le2V+1$, two distinct primes cannot contribute, since the
required powers have total degree at least $2V+2$. The identity
$\max\{t-l_j,0\}^2=\sum_{k=l_j+1}^t(2(k-l_j)-1)$ gives the second moment.
\end{proof}

\subsection{The bilinear sum over all nonzero frequencies}
\begin{lemma}\label{lem:joint-bilinear}
For $r,s\ge1$, $a\in\F_q$ and sequences $A_f,B_g$ on $\M_r,\M_s$,
put
\[
 \mathcal T_a(A,B)=\sum_{t\in\F_q^\times}\psi(-ta)
    \sum_{f\in\M_r,\,g\in\M_s} A_f B_g\psi(t\lam(fg)).
\]
Then
\begin{equation}\label{eq:joint-bilinear}
 |\mathcal T_a(A,B)|\le \beta_{r,s}(a)\|A\|_2\|B\|_2,
 \qquad
 \beta_{r,s}(a)=
 \begin{cases}
  \sqrt{q-1}\,q^{\max(r,s)/2},&r\ne s,\\
  q^{(r+1)/2},&r=s,\ a\ne2,\\
  (q-1)q^{r/2},&r=s,\ a=2.
 \end{cases}
\end{equation}
\end{lemma}
\begin{proof}
By additive orthogonality the kernel is
$K_{f,g}=q\one_{\lam(fg)=a}-1$. Suppose $r<s$ and write $g=Y^s+g_0$.
Distinct monic $f$ in $\M_r$ are nonproportional in the space of
polynomials of degree less than $s$. Lemma~\ref{lem:pairings} implies
that the linear forms $g_0\mapsto\lam(fg_0)$ are nonzero and
pairwise independent. Hence distinct rows of $K$ are orthogonal,
and each row has squared norm $q^s(q-1)$. The case $s<r$ follows
by transposition.

For $r=s$, write $f=D_r+f_0$, $g=D_r+g_0$. Orthogonality in the
Dickson basis gives $\lam(fg)=2+\lam(f_0g_0)$, a translated
dot product as in~\cite[Theorem~2.1]{HIKR}. Nonzero $f_0$ in
distinct one-dimensional subspaces give orthogonal rows. On each such
subspace, the Gram block of the $q-1$ nonzero rows is $q^r(qI-J)$
if $a\ne2$, and $q^r(q-1)J$ if $a=2$, where $J$ is the all-ones
matrix. Indeed, proportional vectors impose distinct nonzero
right-hand sides when $a\ne2$, and the same zero right-hand side
when $a=2$. The row for $f_0=0$ is constant, orthogonal to the other
rows, and has squared norm $q^r$ or $(q-1)^2q^r$, respectively.
Since $J$ has eigenvalues $q-1$ and $0$, the two Gram matrices
have largest eigenvalues $q^{r+1}$ and $(q-1)^2q^r$.
\end{proof}

Let $U,V\ge2$ be integers with
\begin{equation}\label{eq:unit-range}
 2U+3\le n,\qquad U+V+2\le n\le U+2V+2.
\end{equation}
For $r=U+1,\ldots,n-V-1$ set $s=n-r$
and choose $c_r\ge0$. The centred second moment is
\[
 \widetilde M_s=M_s(V)-2c_rS_s(V)+c_r^2V_s.
\]
Let $H(r,m)$ bound the sum of the two masses in
Lemma~\ref{lem:unitmass}, either by its table or by
\eqref{eq:combined-mass}. Define
\begin{align*}
 J_n&=\sum_{\substack{d\mid n\ ,\ n/d>1\text{ odd}}}d\pi_q(d),\\
 \mathcal E_c&=\sum_{r=0}^{U+V}w_rH(r,n-r)+
              \sum_{r=U+1}^{n-V-1}c_rH(r,n-r),\\
 R_c&=\frac{2}{q^n}\sum_{r=U+1}^{n-V-1}
        \beta_{r,n-r}(a)\sqrt{A_r\widetilde M_{n-r}}
     +\frac{q-1}{q^n}\mathcal E_c+2q^{1-n}J_n+q^{-n}.
\end{align*}
Here $w_r=r$, except that one may take $w_{U+V}=V\lfloor(U+V)/V\rfloor$.
\begin{proposition}\label{prop:centred-unitbound}
Under \eqref{eq:unit-range},
\begin{equation}\label{eq:joint-prime}
 N_q(n,a)\ge\frac{q^{n-1}}{2n}(1-R_c).
\end{equation}
\end{proposition}
\begin{proof}
Apply \eqref{eq:Vaughan} to the weighted count
$\tfrac12\sum_{f\in\M_n}\Lambda(f)(\rho(f)-\ch(f))
\one_{\lam(f)=a}$.
The zero frequency is $(q^n-1)/(2q)$ by Lemma~\ref{lem:zero}.
In nonzero frequencies the first Vaughan term vanishes by
Lemma~\ref{lem:unitmass}, since $n\ge2U+3$.
The coefficient of a degree-$r$ product in the second term has
absolute value at most $r$. At $r=U+V$, every contributing
prime power has degree $V$. Distinct such divisors are coprime,
so there are at most $\lfloor r/V\rfloor$, of weight at most $V$
each. This gives the alternative weight $w_{U+V}$.
For each third-term sum, write $B_V(g)=(B_V(g)-c_r)+c_r$.
Here $V<s\le2V+1$ by \eqref{eq:unit-range}. The first summand
has squared norm $\widetilde M_s$ after either unit twist,
so \eqref{eq:joint-bilinear} applies together with the norm
$A_r^{1/2}$ of the twisted M\"obius weight. The second
summand has total mass at most $(q-1)c_rH(r,s)$ over the nonzero
frequencies by Lemma~\ref{lem:unitmass}, since $|\mu(f)|\le1$.
An even power of a prime contributes zero to $\rho-\ch$. Proper
odd powers contribute at most $J_n$ to the weighted count.
Combining these bounds with its zero-frequency contribution proves
\eqref{eq:joint-prime}.
\end{proof}

\subsection{Lower bounds for the remaining ranges}
\begin{proposition}\label{prop:smalldegrees}
For $7\le n\le12$, we have $N_q(n,a)>0$ for all $a\in\F_q$
and all odd prime powers $q\ge q_0$, where
\[
 \begin{array}{c|rrrrrr}
 n&7&8&9&10&11&12\\\hline
 q_0&9&9&5&5&3&3
 \end{array}
\]
The same conclusion holds for $13\le n\le15$ for every odd $q$,
and for $(n,q)=(8,7)$ if $a\ne2$. Over $\F_3$, $N_3(n,a)$ exceeds
one tenth of $3^{n-1}/(2n)$ for $n=16$ and one fifth for $n=17$.
\end{proposition}
\begin{proof}
For the tail bounds take $c_r=0$, $w_r=r$ and the masses bounded
as in Lemma~\ref{lem:unitmass}. We have $A_r\le V_r=q^r(1-q^{-1})^2$
for $r\ge2$. Using $jI_q^*(j)\le q^j$ and $V_{s-kj}\le q^{s-kj}$
in \eqref{eq:moment} gives $M_s(V)\le q^s L_s(V;q)$, where
\[
 L_s(V;q)=\sum_{j=1}^s j\sum_{k=\lfloor V/j\rfloor+1}^{\lfloor s/j\rfloor}
 (2(k-\lfloor V/j\rfloor)-1)q^{-(k-1)j}.
\]
In $R_c$ this bounds each square-root term by $2\sqrt{Z_r}$,
where
\[
 Z_r=\frac{(1-q^{-1})^4q^{2-\min(r,n-r)}L_{n-r}(V;q)}{\delta_r},
 \qquad \delta_r=\begin{cases}q-1&r\ne n-r,\\1&r=n-r.\end{cases}
\]
After removing $L_s^{1/2}$, the factors in $\sqrt{Z_r}$ are
$(1-q^{-1})^{3/2}q^{(1-\min(r,s))/2}$ for $r\ne s$ and
$(1-q^{-1})^2q^{1-r/2}$ for $r=s\ge4$. Both decrease on the
ranges in Table~\ref{tab:bounds}, as does $L_s$, a polynomial in
$q^{-1}$ with nonnegative coefficients.

Write $m=n-r$. After multiplication by
$(q-1)q^{-n}$, each mass bound is $q^{2-m}$ or $q^{1-m}$ times
a factor in $1/q$. For $L=0,1$ this factor is, up to a constant,
$1-q^{-1}$ times one of $1$, $1-2/q$, $1-4/(3q)$.
Its logarithmic derivative with respect to $\log q$ is at most
$1/(q-1)+2/(q-2)<1$ for $q\ge5$, and $m\ge4$. For $L=2,3$ or
$L\ge4$ the factor is $(1-q^{-1})^k$ with $k=3$ or $5$,
respectively. On the listed ranges $k/(q-1)<m-2$ (the case
$L\ge4$, $m=3$ has $q\ge17$). Thus each normalized mass decreases
with $q$. The terms $q^{1-n+d}$ bounding the contribution of
proper prime powers decrease as well.

For a single field we use \eqref{eq:joint-prime} directly. In
\eqref{eq:combined-mass}, any consistent intersection of the four
affine conditions $x=0,y=0,u=0,v=0$ has $q^{L-j}$ elements if
its rank is $j$. For example, at $(n,q,r,m)=(7,9,4,3)$ write
$C=c_1T+c_2T^2+c_3T^3+T^4$. Then $u=2+c_3-c_1$ and
$v=2-c_3+c_1$, while for fixed $c_3-c_1$ the pair $(x,y)$ runs
freely through $\F_9^2$. Formula~\eqref{eq:combined-mass} gives
$(q-2)(q-1)^2(q^2+4q)=52416$. The other exact masses used below
are included in the following list, in increasing order of
$r$ starting at $\lceil(n-2)/2\rceil$.
\[
\begin{array}{rrl}
\toprule
n&q&H(r,n-r)\\\midrule
7&9&(6966,52416)\\
8&7&(2744,14504,99792)\\
9&5&(2500,11450,57600)\\
11&3&(324,972,3024)\\
12&3&(0,1620,3024,9072)\\\bottomrule
\end{array}
\]

In Table~\ref{tab:bounds}, a row marked $\dagger$ applies to the
listed $q$ and all larger odd prime powers, using the majorant $Z_r$
above and $d\pi_q(d)\le q^d$ in $J_n$. In the single-field rows take
$Z_r=q^{-2n}\beta_{r,n-r}(a)^2 A_r\widetilde M_{n-r}$.
Rows marked $e$ use exact masses and the alternative $w_{U+V}$.
The centres are $c_r=n-r-V$ in rows $*$ and $e$, except that
$(n,q)=(7,9),(11,3)$ use $(2,1/3),(4,5/2,3/2)$, respectively.
Rows marked $0$ use $c_r=0$. Only the row $(8,7)$ requires $a\ne2$.
In each row, in increasing order of $r$,
\[
 Z_r<(k_r/10^4)^2,\qquad
 (q-1)q^{-n}\mathcal E_c+2q^{1-n}J_n+q^{-n}\le t/10^4 .
\]
\begin{center}
\begin{minipage}{\linewidth}
\centering
\setlength{\abovecaptionskip}{0pt}
\captionof{table}{Upper bounds for the terms in \eqref{eq:joint-prime}.}\label{tab:bounds}
\smallskip
\small
\begin{tabular}{crrrrlr}
\toprule
&$n$&$q$&$U$&$V$&$(k_r)$&$t$\\\midrule
$\dagger$&7&13&2&2&$(1808,1183)$&3338\\
$\dagger$&8&17&2&3&$(1612,1044)$&2832\\
$\dagger$&9&9&3&2&$(1078,824,1617)$&679\\
$\dagger$&10&7&3&3&$(1663,1194,862)$&1051\\
$\dagger$&11&5&4&3&$(1114,864,1295)$&2161\\
$\dagger$&12&5&4&4&$(1216,851,641)$&2322\\
$\dagger$&13&5&5&4&$(544,426,641)$&2621\\\midrule
$e$&7&9&2&2&$(1365,1376)$&4382\\
$*$&7&11&2&2&$(1191,1020)$&5055\\
$e$&8&7&2&2&$(1501,644,1398)$&2116\\
$*$&8&9&2&2&$(1309,1284,1183)$&1724\\
$*$&8&11&2&2&$(1151,1134,1020)$&1323\\
$*$&8&13&2&2&$(1024,1011,895)$&1067\\
$e$&9&5&3&2&$(753,817,1674)$&2952\\
$*$&9&7&3&2&$(565,597,1397)$&2640\\
$*$&10&5&3&3&$(1161,1001,953)$&2276\\
$e$&11&3&4&3&$(834,854,1290)$&3824\\
$e$&12&3&4&3&$(846,731,915,1430)$&1809\\
$*$&13&3&4&4&$(1185,678,640,1025)$&2013\\
$*$&14&3&5&4&$(675,557,638,1026)$&1924\\
$*$&15&3&5&5&$(827,463,430,665)$&1885\\\midrule
$0$&16&3&6&5&$(773,543,554,776)$&1850\\
$0$&17&3&6&6&$(822,411,343,490)$&2017\\\bottomrule
\end{tabular}
\end{minipage}
\end{center}
In each row $2\sum_r k_r+t<10^4$, so $R_c<1$. The last two rows give
$R_c<7142/10^4$ and $R_c<6149/10^4$, respectively. For $n=14,15$
and $q\ge5$ use \eqref{eq:large-ranges}.
\end{proof}

\section{Reduced degrees five and six}\label{sec:fivesix}

For $b\in\M_r$, $m\ge1$, let
\[
 \Delta_\theta(b;m)=\sum_{h\in\M_m}\theta(bh)
                    \bigl(\one_{\lam(bh)=a}-q^{-1}\bigr).
\]
As in Lemma~\ref{lem:mass}, $\Delta_1$ vanishes unless
$\lam(bY^j)=0$ for $0\le j<m$. When $r<m$ no monic $b$ qualifies.
For $r=m$ only $D_m$ qualifies, and the absolute deviation is at most
$q^{m-1}(q-1)$. When $r>m$, nondegeneracy on $\F_q[Y]_{<m+1}$
makes the $m+1$ values $\lam(bY^j)$, $0\le j\le m$, independent
affine coordinates. Among the $q^{r-m}$ surviving choices of $b$,
exactly $q^{r-m-1}$ satisfy $\lam(bY^m)=a$. Thus
$\sum_{b\in\M_r}|\Delta_1(b;m)|=2q^{r-2}(q-1)$.

For $\theta=\ch$ and $m\ge2$, use the coordinates of \eqref{eq:evaluation-line}.
The unrestricted sum is zero because the endpoint values $u,v$
are independent. Averaging over $k\in\M_{m-2}$ also kills the
deviation unless $\lam(Dbw)=0$ for every $\deg w<m-2$. In that case
$\lam(Dbk)=\lam(DbY^{m-2})$ and
\begin{equation}\label{eq:lifted-line}
 \Delta_\ch(b;m)=q^{m-2}\ch(b)
       \sum_{\alpha_bu+\beta_bv=a-\lam(DbY^{m-2})}\eta(uv).
\end{equation}
If $\alpha_b\beta_b\ne0$, the line sum equals $-\eta(-\alpha_b\beta_b)$
when its right-hand side is nonzero and $(q-1)\eta(-\alpha_b\beta_b)$
otherwise. This is the quadratic character sum in $u$ after eliminating
$v$. When $\alpha_b\beta_b=0$ the sum is zero.

\subsection{Reduced degree five}
\begin{proposition}\label{thm:five}
For odd $q$, $D_q(5)\ll q^{-1/2}$. If $q\ge11$, then
$N_q(5,a)>0$ for every $a\in\F_q$.
\end{proposition}
\begin{proof}
Apply \eqref{eq:Vaughan} with $(U,V)=(1,2)$ to the deviation of
$W_5(a)$ from its mean. For $\deg b\le1$ both deviations vanish
by Lemma~\ref{lem:mass}, eliminating the first term and these
degrees in the second.

In the second convolution write
\[
 \gamma(b)=\sum_{de=b,\,\deg d\le1,\,\deg e\le2}\mu(d)\Lambda(e).
\]
For quadratics this equals $-2\mu(b)$, by considering irreducibles,
split squarefree polynomials, and squares. The untwisted deviation
vanishes for a quadratic $b$ by Lemma~\ref{lem:pairings}.
For a cubic $b$ the factor $e$ must be a quadratic prime power.
Two distinct such divisors would be coprime and have total degree
four, so $|\gamma(b)|\le2$. The bound for $r=3$, $m=2$ gives untwisted
error at most $4q(q-1)$.

For the twisted cubic term sum over the contributing pairs
$b=(Y-r)(Y^2+v_1Y+z)$, with weight bounded by
$\Lambda(Y^2+v_1Y+z)$. Only $r\ne\pm2$
can contribute, giving total weight at most $(q-2)q^2$.
Moreover the equation $\lam(Db)=a$ gives $v_1=r+rz-a/2$.
For fixed $r\ne0$, the discriminant $v_1^2-4z$ is quadratic in
$z$ with nonzero square leading coefficient, so the sum of
$1-\eta(v_1^2-4z)$ over $z$ is at most $q+1$. For $r=0$ that
sum is $q$. Since $\Lambda(Y^2+v_1Y+z)=1-\eta(v_1^2-4z)$,
the weight on $\lam(Db)=a$ is at most $q+(q-3)(q+1)=q^2-q-3$.
The line sum in \eqref{eq:lifted-line} therefore bounds the twisted
cubic contribution by
\[
 (q-2)q^2+(q-2)(q^2-q-3)=2q^3-5q^2-q+6 .
\]

The third convolution has degrees $2+3$. Write $B_2(g)$ as
$(B_2(g)-1)+1$ on cubics in each twist. The constant part
of the untwisted deviation vanishes by the pairing.
Combining its twisted counterpart with the quadratic contribution
of the second convolution gives weight $\mu(b)-\gamma(b)=3\mu(b)$.
For $b=Y^2+sY+t$, \eqref{eq:lifted-line} with $m=3$ requires
$\lam(Db)=0$, or $t=-1$. Then $\lam(DbY)=-2s$, so only one
choice of $s$ gives the larger line sum. The total error is at most
$3q((q-1)+(q-1))=6q(q-1)$.

For the centred term, $B_2(g)=\Lambda(g)$ on cubics. It equals $3$
on irreducible cubics, $1$ on cubes of linear polynomials and $0$
otherwise. Thus
\begin{equation}\label{eq:cubic-moment}
 \sum_{g\in\M_3}(B_2(g)-1)^2=9\pi_q(3)+q-q^3=2q^3-2q .
\end{equation}
The squared norm of $\mu(b)$ on quadratics is $q^2-q$. Both norms can only decrease
under the twist. Dividing the unequal-degree estimate
\eqref{eq:joint-bilinear} by $q$ gives the operator bound
$\sqrt{q(q-1)}$ for the deviation kernel. Thus the two centred
contributions together have error at most
$2\sqrt2(1-q^{-1})q^{7/2}$.

The only proper powers are fifth powers of linears, of total weight
at most $2q$ in $W_5$. The cubic and quadratic errors above, together
with these powers, contribute at most $2q^3+5q^2-9q+6$.
The mean of $W_5$ is $q^4+q^{-1}$. Dividing by $q^4$ gives
\[
 D_q(5)\le\frac{2\sqrt2(1-1/q)}{\sqrt q}+\frac{2}{q}+\frac5{q^2}.
\]
The right side decreases for real $q\ge11$, and at $q=11$ it is
less than $44865/44891$, using $\sqrt2<99/70$ and
$\sqrt{11}>53/16$.
\end{proof}

\subsection{Reduced degree six}
\begin{proposition}\label{prop:six}
For odd $q$, $D_q(6)\ll q^{-1/2}$. One has $N_q(6,a)>0$ if
$q\ge13$, or if $q\ge9$ and $a\ne2$.
\end{proposition}
\begin{proof}
Use \eqref{eq:Vaughan} with $U=V=2$ in $W_6(a)$.
In the untwisted first two convolutions, factors of degree
$r\le2$ have zero deviation. The coefficient of $b$ in the
second convolution has absolute value at most $r$. The untwisted
bounds for $(r,m)=(3,3),(4,2)$ therefore give errors at most
$3q^2(q-1)$ and $8q^2(q-1)$, respectively.

For the twist $\ch$, the $D$-pairing and \eqref{eq:evaluation-line}
give zero deviation for $r\le1$. When $r=2$, only $b=U_2=Y^2-1$
can contribute by the $D$-pairing, with combined weight at most
$4+2$ in the first two convolutions. Its error is bounded by
$6q^2(q-1)$ using \eqref{eq:lifted-line}. For $r=3$, the
condition $\lam(Db)=0$ restricts $b$ to $Y^3+sY^2+tY-s$.
Since $\lam(YDb)=-4-2t$, the right-hand side of the line constraint
in \eqref{eq:lifted-line} is zero for exactly $q$ of these $q^2$
cubics. The bound is therefore
$3q[q(q-1)+(q^2-q)]=6q^2(q-1)$.

For $r=4$ write $b=de$, where $d=Y^2+sY+t$ and $e=Y^2+vY+z$.
The right side of the constraint on the evaluations of the
remaining quadratic is $a-\lam(Dde)$, with
\[
 \lam(Dde)=-4-2(t+z+sv)-2tz.
\]
The total weight $|\mu(d)|\Lambda(e)$ is at most $q^4$ by
\eqref{eq:zerofreq}. For $t\ne-1$ and fixed $(s,t)$, the condition
$\lam(Dde)=a$ specifies $z$ as an affine function of $v$. Hence
$v^2-4z$ is monic quadratic in $v$ and $\sum_v\Lambda(e)\le q+1$.
For $t=-1$ the condition reads $-2-2sv=a$. If $s\ne0$,
$v$ is fixed and summing $1-\eta(v^2-4z)$ over $z$ gives $q$.
If $s=0$, all $e$ together have weight $q^2$. Thus the exceptional
locus has weight at most $q(q-1)(q+1)+q(q-1)+q^2=q^3+2q^2-2q$.
By \eqref{eq:lifted-line} its contribution and that of its complement
are bounded in total by
\[
 q^4+(q-2)(q^3+2q^2-2q)\le2q^4.
\]
The combined first and second convolution errors are bounded
by $2q^4+23q^2(q-1)$.

In the sole split $3+3$ of the third convolution, write
$B_2=(B_2-1)+1$. The constant part contributes at most $3q^2(q-1)$
for the two twists, by the degree-three fibre bounds above with
coefficient one. For the centred part use \eqref{eq:cubic-moment},
$\sum_{d\in\M_3}\mu(d)^2=q^3-q^2$, and Lemma~\ref{lem:bilinear}.
Summing the $q-1$ nonzero frequencies and dividing by $q$ bounds
the two twists together by $2\sqrt2(1-q^{-1})^{3/2}q^{9/2}$.
Proper prime powers contribute to the difference of the twists
at most $2(q^2-q+1)$. Indeed, the cubes of irreducible quadratics
contribute at most $2(q^2-q)$, and even powers contribute only
from $(Y\pm2)^6$. The mean of $W_6(a)$ is $q^5+q^{-1}$. Since
$2(q^2-q+1)+q^{-1}\le 2q^2$, the absolute errors above give
\[
 D_q(6)\le\frac{2\sqrt2(1-q^{-1})^{3/2}}{\sqrt q}
 +\frac2q+\frac{26(1-q^{-1})}{q^2}+\frac2{q^3}.
\]
The right side decreases for $q\ge13$. At $q=13$, the inequalities
$\sqrt2<17/12$, $\sqrt{13}>18/5$, and $\sqrt{1-1/q}\le1-1/(2q)$
give an upper bound of $996/1000$.

For $a\ne2$, use the uncentred $B_2(g)$ with squared norm at most
$3q^3$. The balanced case of \eqref{eq:joint-bilinear} bounds the
third-convolution error by $2\sqrt3q^4$ after dividing the sum
over frequencies by $q$. The sufficient inequality becomes
\[
 \frac{2+2\sqrt3}{q}+\frac{23}{q^2}+\frac2{q^3}<1.
\]
At $q=9$ the left side is less than $3259/3645$ by $\sqrt3<26/15$,
and it decreases thereafter.
\end{proof}

\begin{proof}[Proof of Theorem~\ref{thm:distribution}]
For fixed $q$, \eqref{eq:distributionexplicit} gives the additive
error $O_q(q^{3n/4})$. For $5\le n\le9$, Propositions~\ref{thm:five}
and~\ref{prop:six} and \eqref{eq:relative-seven-nine} give
$D_q(n)\ll q^{-1/2}$. For $n\ge10$, put $m=\lfloor n/2\rfloor$.
Since $m\ge5$ and $C_q<q^2$, \eqref{eq:relative} gives
\[
 D_q(n)\ll (2m+1)q^{2-m/2}
 \le \min\bigl\{(2m+1)3^{2-m/2},\,11q^{-1/2}\bigr\}.
\]
Indeed, for $m\ge5$ the ratio of successive terms in $m$ is
$(2m+3)/((2m+1)\sqrt q)<1$. The first bound tends uniformly to
zero with $n$, the second with $q$, proving the joint limit.
For the explicit inequalities, \eqref{eq:large-ranges} covers
$q\ge5$, $n\ge16$, and $q=3$, $n\ge18$.
Proposition~\ref{prop:smalldegrees} covers the remaining two pairs.
\end{proof}

\section{Artin--Schreier constructions in degree \texorpdfstring{$2p$}{2p}}\label{sec:AS}
Over any finite field $\F_q$ of characteristic $p$, write
$\Tr=\Tr_{\F_q/\F_p}$. A root $\zeta$ of $Z^p-Z+t$ satisfies
$\zeta^q=\zeta-\Tr(t)$. Frobenius permutes the $p$ roots by a
translation, so
\begin{equation}\label{eq:AScriterion}
 Z^p-Z+t\text{ is irreducible over }\F_q
 \quad\Longleftrightarrow\quad \Tr(t)\ne0 .
\end{equation}
This criterion for translation-invariant polynomials
(cf.~\cite[Section~2]{KapetanakisInvariant}) supplies the extension
of reduced degree $p$.

\begin{theorem}\label{thm:AS}
Let $q=p^e$ with $p\ge5$ prime. Then $N_q(p,a)>0$ for every $a\in\F_q$.
If $q=p$, choose $v\in\F_p$, $v\ne a$, with $\eta_p(v^2-4)=-1$,
where $\eta_p$ is the quadratic character of $\F_p$.
An explicit polynomial with central coefficient $a$ is
\begin{equation}\label{eq:ASprime}
 B_v^p+\frac{a-v}{2}\bigl(X^p-XB_v^{p-1}\bigr),
 \qquad B_v=X^2+vX+1 .
\end{equation}
\end{theorem}
\begin{proof}
Take $w,t\in\F_q^\times$ with $\Tr(t)\ne0$, and let $\zeta$ be
a root of $Z^p-Z+t$. For $v\in\F_p\setminus\{\pm2\}$, the element
$w/\zeta-v$ generates the same extension as $\zeta$ and has minimal polynomial
\[
 G(Y)=(Y+v)^p-\frac wt(Y+v)^{p-1}+\frac{w^p}{t}.
\]
Its evaluations are $G(\pm2)=v\pm2+(w^p-w)/t$. For the central
coefficient we use
\begin{equation}\label{eq:centralp}
 [X^{p-1}]B_v^{p-1}=(v^2-4)^{(p-1)/2}.
\end{equation}
Indeed, the coefficient of $v^{p-1-2j}$ on the left is congruent
modulo $p$ to $\binom{2j}{j}=(-4)^j\binom{(p-1)/2}{j}$ in $\F_p$,
the corresponding coefficient on the right.

When $q=p$, there are $(p-1)/2\ge2$ choices of $v$ with
$\eta_p(v^2-4)=-1$, so one can exclude $a$. Set $w=1$ and
$t=2/(a-v)$. Now $G(2)G(-2)=v^2-4$ is a nonsquare, and
Lemma~\ref{lem:transform} gives \eqref{eq:ASprime}, with central
coefficient $v+(a-v)=a$ by \eqref{eq:centralp}.

For the extension fields we need the following consequence of
\eqref{eq:Weileven}. If $q\ge p^2$, $A,B\in\F_q^\times$ and
$\delta=\pm1$, then
\begin{equation}\label{eq:trace-select}
 \Tr(u)\ne0,\qquad \eta(A+Bu^{p-1})=\delta
\end{equation}
has a solution, provided $\eta(A)\ne\delta$ when $q=p^2$.
Indeed, $A+BU^{p-1}$ is squarefree, has at most $p-1$ zeros, and
its character sum has absolute value at most $(p-3)\sqrt q+1$.
After omitting the $q/p$ trace-zero elements, at least
\[
 \tfrac12\bigl((1-2/p)q-p-(p-3)\sqrt q\bigr)
 =\tfrac12(\sqrt q-p)\bigl((1-2/p)\sqrt q+1\bigr)
\]
choices remain. At $q=p^2$, $u=0$ was already excluded by its
character, improving this lower bound by one.

Suppose $q\ge p^2$ and $a\ne\pm2$. Put
$k=(2-a)/(a+2)$, $A=1+k^{1/p}$ and $\kappa=k/A^{p-1}$.
Since $k\ne0,-1$, both $A$ and $\kappa$ are nonzero.
For $s=A/u$, $d=s^{p-1}+k$, set $J=T^p-s^{p-1}T+d$. Then
\[
 \Tr(d/s^p)=\Tr((1+k^{1/p})/s)=\Tr(u),\qquad
 \eta(d)=\eta(1+\kappa u^{p-1}).
\]
Use \eqref{eq:trace-select} to obtain $\Tr(u)\ne0$ and
$\eta(d)=-\eps$. At $q=p^2$ the desired sign is $-1$, distinct
from $\eta(1)$. Equation~\eqref{eq:AScriterion} applied to
$s^{-p}J(sZ)$ gives irreducibility of $J$, and $\eta(-d)=-1$
allows us to apply Lemma~\ref{lem:cayley}.
The coefficients of $X^p$ in $(X+1)^{2p}$, $(X-1)^{2p}$ and
$(X+1)^2(X-1)^{2p-2}$ are $2,-2,-2$ in characteristic $p$.
For the last identity, write
$(X+1)^2(X-1)^{2p-2}=(X-1)^{2p}+4X(X-1)^{2p-2}$ and use
$\binom{2p-2}{p-1}\equiv0\pmod p$.
Thus the transformed central coefficient is
\[
 \frac{2+2s^{p-1}-2d}{1-s^{p-1}+d}
 =\frac{2-2k}{1+k}=a.
\]

By \eqref{eq:symmetry} only $a=2$ remains. If $\eps=-1$, take
$k=0$, $\Tr(1/s)\ne0$ and $d=s^{p-1}$. The same construction
gives $a=2$ since $\Tr(d/s^p)\ne0$ and $\eta(-d)=-1$.
If $\eps=1$, choose nonzero $v\in\F_p$ with $\eta_p(v^2-4)=-1$.
By \eqref{eq:trace-select} there is $w$ of nonzero trace with
\begin{equation}\label{eq:AScondition}
 \eta\left(-\frac4{2-v}(2w^{p-1}+v)\right)=-1.
\end{equation}
At $q=p^2$ the nonzero constant $-4v/(2-v)\in\F_p$ is a square,
as required. Set $t=(w^p+w)/(2-v)$, whose trace is
$2\Tr(w)/(2-v)$. In the polynomial $G$ at the start of the proof,
these parameters give
\[
 G(2)G(-2)=-\frac{4w^2(2w^{p-1}+v)}{t^2(2-v)}.
\]
By \eqref{eq:AScondition}, this is a nonsquare. The transform of
$G$ is $B_v^p-(w/t)XB_v^{p-1}+(w^p/t)X^p$, whose central
coefficient is $v+(w^p+w)/t=2$ by \eqref{eq:centralp}.
\end{proof}

\section{Completion of the existence proof}\label{sec:finite}
\begin{proof}[Proof of Theorem~\ref{thm:classification}]
The quadratic discriminant and Propositions~\ref{prop:n2},
\ref{prop:n3} and~\ref{prop:n4} settle $n\le4$. For $n\ge5$, apply
Theorem~\ref{thm:distribution} and Propositions~\ref{prop:smalldegrees},
\ref{thm:five} and~\ref{prop:six}. Theorem~\ref{thm:AS} covers
$(q,n)=(5,5),(7,7)$. The remaining cases are
\[
\begin{array}{c|l}
n&\text{Small-field cases}\\\hline
5&q=3,7,9\\
6&q=3,5,7\ ;\ q=9,11,\ a=2\\
7&q=3,5\\
8&q=3,5\ ;\ q=7,\ a=2\\
9,10&q=3
\end{array}
\]

We use the monic transforms $H[\tau]$ defined before
Lemma~\ref{lem:cayley}. For each resulting irreducible $G$ it
suffices to verify $\lam(G)=a$ and that $\Pi=G(2)G(-2)$ is a
nonsquare, by Lemma~\ref{lem:transform}.

\smallskip\noindent\emph{The field $\F_3$.}
For $5\le n\le10$, use the base polynomials
\begin{align*}
 H_5&=Y^5+Y^4-Y^3+Y^2-1, & H_7&=Y^7+Y^6+Y^5+Y^4+1,\\
 K_3&=Y^3-Y^2+Y+1, & L&=Y^4+Y^3-1,\\
 M&=Y^4+Y^2-1.
\end{align*}
The product of $H_5$ and its monic reciprocal is $\Phi_{11}$, whose
irreducible factors have degree $\operatorname{ord}_{11}(3)=5$.
The cubic $K_3$ and the quartics $L,M$ have no rational roots.
The only products of irreducible quadratics with constant term $-1$
are $(Y^2+1)(Y^2\pm Y-1)$, neither of which equals $L$ or $M$.

The polynomial $H_7$ has no rational root. Modulo $H_7$,
the residues of $Y^9-Y$ and $Y^{27}-Y$ are respectively
$Y^2(Y-1)^3$ and $Y^3(Y^2-Y-1)$. The remainder of $H_7$ modulo
$Y^2-Y-1$ is $1-Y$. Thus neither residue has a common factor with
$H_7$. This excludes all factors of degree at most three and
proves irreducibility.

For the sextics, let $\delta^2=-1$ and take
$K=Y^3+Y^2+tY+u$ for $(t,u)=(1,-1),(-1,1)$.
Neither cubic has a root in $\F_3$, and hence both remain irreducible
over $\F_9$. The two cubics $K(Y\pm\delta)$ are distinct, since
their coefficients of $Y$ differ by $4\delta$. Being conjugate,
they have irreducible product over $\F_3$. For the two choices
of $(t,u)$ the products are respectively
\[
 G_2=Y^6-Y^5+Y^3+Y^2-Y+1,\qquad
 G_1=Y^6-Y^5-Y^4+Y^3-Y^2+Y+1.
\]
Their transforms have central coefficients $\lam(G_2)=2$ and
$\lam(G_1)=1$. The monic reciprocal $G_2^*$ of $G_2$ gives
$\lam(G_2^*)=0$. In all three cases $G(1)G(-1)=-1$.

Let $K_5=H_5[1/(Y-1)]$, so $K_5(0)=1$. The square-root criterion
in Lemma~\ref{lem:cayley} gives the irreducibility of $K_5(Y^2)$,
$L(Y^2)$ and $M(Y^2)$. A root $u$ of $K_3$ has trace $1$ to
$\F_3$. Thus $Y^3-Y-u$ is irreducible over $\F_{27}$ by
\eqref{eq:AScriterion}, and $K_3(Y^3-Y)$ is irreducible of degree nine.
For each displayed substitution $\tau$, $G=H[\tau]$ has $\Pi=-1$.
Writing $g_j=[Y^j]G$, \eqref{eq:lambda} reduces here to
$\lam(G)=g_0-g_2-g_6+g_8$. In odd degree the sign symmetry
\eqref{eq:symmetry} supplies the opposite value.
\[
\begin{array}{c|l|ccc}
 n&H&a=0&a=1&a=2\\\hline
 5&H_5&1/(Y+1)&1/Y&\\[1pt]
 7&H_7&1/(Y+1)&-1/(Y+1)&\\[1pt]
 8&L(Y^2)&(Y+1)/(Y-1)&&Y+1\\[1pt]
 8&M(Y^2)&&(Y+1)/(Y-1)&\\[1pt]
 9&K_3(Y^3-Y)&1/Y&1/(Y+1)&\\[1pt]
 10&K_5(Y^2)&Y/(Y+1)&(Y+1)/Y&1/(Y+1)\\[1pt]
\end{array}
\]
\smallskip\noindent\emph{The other fields.}
For a prime $\ell=2m+1$ not dividing $q$, the polynomial
$1+\sum_{j=1}^{m}D_j$ has roots $z+z^{-1}$ for primitive
$\ell$th roots $z$, since evaluation gives $z^{-m}\Phi_\ell(z)$.
The degree of each such root over $\F_q$ is the least $d\ge1$
with $q^d\equiv\pm1\pmod\ell$, since $z+z^{-1}$ determines
the unordered pair $\{z,z^{-1}\}$.

\smallskip\noindent\emph{Reduced degree five.}
Over $\F_7$, take $H=1+\sum_{j=1}^5D_j$. Its roots have
degree five since $7^5\equiv-1\pmod {11}$ and $7\not\equiv\pm1\pmod {11}$.
Over $\F_9=\F_3(i)$, where $i^2=-1$, take $H=Y^5+Y^4+Y^2+1$.
The latter has no root in $\F_9$. Indeed, for $y^4=-1$ it takes
the value $y(y-1)$, and for $y\in\{\pm1,\pm i\}$ it equals
$y^2+y-1\ne0$. Also $H(0)=1$. Thus it has no factor of degree
one or two over $\F_3$, and is irreducible of degree five over
both fields. For $G=s^5H((Y-t)/s)$ the entries below give all
values of $a$ up to sign and Frobenius conjugation.
\[
\begin{array}{cc|cc@{\qquad\qquad}cc|cc}
\multicolumn{4}{c}{q=7}&\multicolumn{4}{c}{q=9}\\
s&t&a&\Pi&s&t&a&\Pi\\\hline
1&3&0&3&1+i&-1&0&-1+i\\
1&0&1&3&1+i&0&1&-1+i\\
3&2&5&3&1&i&i&-1-i\\
3&4&3&6&1+i&-i&1+i&1-i
\end{array}
\]
The entries $\Pi$ are nonsquares. For $q=9$ each has norm $2$ to $\F_3$.

\smallskip\noindent\emph{Reduced degree six.} Take
\[
\begin{array}{c|l}
q&H_q(Y)\\\hline
5&1+Y+Y^2+Y^3+Y^4+Y^5+Y^6\\
7&(Y^2+5)^3-2\\
9&(Y^3-Y)^2-(Y^3-Y)-i,\qquad i^2=-1\\
11&1+\displaystyle\sum_{j=1}^6D_j(Y)\\
\end{array}
\]
For $q=5,11$ irreducibility follows from $\operatorname{ord}_7(5)=6$
and $\min\{d\ge1:11^d\equiv\pm1\pmod{13}\}=6$, respectively.
For $q=7$, a cube root $u$ of the noncube $2$ has degree three,
and $u-5$ has nonsquare norm $2-5^3=3$ in $\F_7$. A square root
of $u-5$ therefore has degree six. Over $\F_9$, $T^2-T-i$ is
irreducible since its discriminant $1+i$ has norm $2$ to $\F_3$.
A root $u$ has trace $1$ to $\F_9$ and $2$ to $\F_3$, so
$Y^3-Y-u$ is irreducible over $\F_{81}$ by \eqref{eq:AScriterion}.

For $G=H_q[\varphi]$ we have $\lam(G)=20+6g_4+2g_2+g_0$.
The following substitutions realize the indicated values of $a$.
\[
\begin{array}{c|lc|lc}
 &\multicolumn{2}{c|}{q=5}&\multicolumn{2}{c}{q=7}\\
 a&\varphi(Y)&\Pi&\varphi(Y)&\Pi\\\hline
 0&2Y&2&1/(Y+2)&3\\
 1&2Y+4&3&Y+3&3\\
 2&Y+3&2&Y/(Y+2)&6\\
 3&Y+1&3&1/(Y+1)&5\\
 4&1+1/Y&2&Y+1&6\\
 5&&&2+1/Y&3\\
 6&&&Y+2&6
\end{array}
\]
For $a=2$ over $\F_9$ take $\varphi=(1+i)Y$. Then
$(g_4,g_2,g_0)=(i,2,2)$ and $\Pi=2+i$ has norm $2$ to $\F_3$.
Over $\F_{11}$, $\varphi=Y+7$ gives $(6,2,8)$ and $\Pi=10$.

\smallskip\noindent\emph{Reduced degrees seven and eight.}
The three remaining pairs $(q,n)$ use the base polynomials
\[
\begin{array}{ccl}
q&n&H(Y)\\ \hline
5&7&Y^7+2Y^6+2Y^5-2Y^4-2Y^3+2Y^2-Y-1\\
5&8&Y^8+2Y^6+Y^4+2\\
7&8&Y^8+Y^4-1 .
\end{array}
\]
For the septic $H$ over $\F_5$, the product of $H$ with
$Y^7-Y^6+2Y^5-2Y^4+Y^3+2Y^2+2Y+1$ is
$1+\sum_{j=1}^{14}D_j$, by the recurrence for the Dickson polynomials.
The roots of this product have degree seven, since
$5^7\equiv-1\pmod{29}$ and $5\not\equiv\pm1\pmod{29}$. Thus $H$ is irreducible.

For the octic over $\F_5$, write $H(Y)=L(Y^2)$, where
$L(Z)=(Z^2+Z)^2-3$. If $u^2=3$, the discriminant $1+4u$ of
$Z^2+Z-u$ has nonsquare norm $3$ in $\F_5$. Hence a root $v$ of
$Z^2+Z-u$ has degree four and minimal polynomial $L$. Its norm
equals $L(0)=2$, so a square root of $v$ has degree eight.
Over $\F_7$, the four polynomials $T^2+bT-1$ with $b=\pm1,\pm3$
have product $T^8+1$. Substituting $T=Y^4$ shows that $H$ divides
$Y^{32}+1$. Every root of $Y^{32}+1$ has order $64$, and
$\operatorname{ord}_{64}(7)=8$, proving the claim.

Over $\F_5$ the substitutions below give $G=H[\varphi]$ with
$\lam(G)=a$. In degree seven the sign symmetry supplies $a=1,3$.
\[
\begin{array}{c|lc|lc}
 &\multicolumn{2}{c|}{n=7}&\multicolumn{2}{c}{n=8}\\
 a&\varphi(Y)&\Pi&\varphi(Y)&\Pi\\\hline
 0&(Y+1)/(2Y+1)&3&1/(Y+1)&3\\
 1&&&Y/(Y+2)&2\\
 2&1/Y&3&1/(2Y+1)&3\\
 3&&&(Y+1)/(Y+2)&2\\
 4&1/(Y+4)&2&1/(2Y+2)&3
\end{array}
\]
Finally, for $(q,n)=(7,8)$ take $\varphi=1/(Y+1)$.
Here $(g_6,g_4,g_2,g_0)=(0,6,1,6)$, giving $\lam(G)=2$ and
$\Pi=3$.
\end{proof}

%


\begin{thebibliography}{99}
\bibitem{BassaMenares}
A. Bassa and R. Menares, The $R$-transform as power map and its
generalizations to higher degree,
\emph{Finite Fields Appl.} \textbf{102} (2025), Article 102546.
\href{https://doi.org/10.1016/j.ffa.2024.102546}{doi:10.1016/j.ffa.2024.102546}.
\bibitem{CarlitzQuartic}
L. Carlitz, A special quartic congruence, \emph{Math. Scand.}
\textbf{4} (1956), 243--246.
\href{https://doi.org/10.7146/math.scand.a-10472}{doi:10.7146/math.scand.a-10472}.

\bibitem{Carlitz}
L. Carlitz, Some theorems on irreducible reciprocal polynomials over a finite
field, \emph{J. Reine Angew. Math.} \textbf{227} (1967), 212--220.
\href{https://doi.org/10.1515/crll.1967.227.212}{doi:10.1515/crll.1967.227.212}.
\bibitem{CiglerTyson}
J. Cigler and M. Tyson, An interesting class of Hankel determinants,
in V. Pillwein and C. Schneider (eds.),
\emph{Algorithmic Combinatorics: Enumerative Combinatorics, Special
Functions and Computer Algebra}, Springer, 2020, pp.~137--157.
\href{https://doi.org/10.1007/978-3-030-44559-1_8}{doi:10.1007/978-3-030-44559-1\_8}.

\bibitem{CohenR}
S. D. Cohen, The explicit construction of irreducible polynomials over
finite fields, \emph{Des. Codes Cryptogr.} \textbf{2} (1992), 169--174.
\href{https://doi.org/10.1007/BF00124895}{doi:10.1007/BF00124895}.
\bibitem{Gao22}
Z. Gao, Enumeration of self-reciprocal irreducible monic polynomials with
prescribed leading coefficients over a finite field,
\emph{Finite Fields Appl.} \textbf{83} (2022), Article 102083.
\href{https://doi.org/10.1016/j.ffa.2022.102083}{doi:10.1016/j.ffa.2022.102083}.
\bibitem{Gao23}
Z. Gao, New estimates and existence results about irreducible polynomials
and self-reciprocal irreducible polynomials with prescribed coefficients
over a finite field, \emph{La Matematica} \textbf{2} (2023), 789--815.
\href{https://doi.org/10.1007/s44007-023-00062-1}{doi:10.1007/s44007-023-00062-1}.
\bibitem{Garefalakis11}
T. Garefalakis, Self-reciprocal irreducible polynomials with prescribed
coefficients, \emph{Finite Fields Appl.} \textbf{17} (2011), 183--193.
\href{https://doi.org/10.1016/j.ffa.2010.11.004}{doi:10.1016/j.ffa.2010.11.004}.
\bibitem{GK12}
T. Garefalakis and G. Kapetanakis, On the Hansen--Mullen conjecture for
self-reciprocal irreducible polynomials,
\emph{Finite Fields Appl.} \textbf{18} (2012), 832--841.
\href{https://doi.org/10.1016/j.ffa.2012.03.001}{doi:10.1016/j.ffa.2012.03.001}.
\bibitem{GK15}
T. Garefalakis and G. Kapetanakis, A note on the Hansen--Mullen conjecture
for self-reciprocal irreducible polynomials,
\emph{Finite Fields Appl.} \textbf{35} (2015), 61--63.
\href{https://doi.org/10.1016/j.ffa.2015.03.005}{doi:10.1016/j.ffa.2015.03.005}.
\bibitem{HamMullen}
K. H. Ham and G. L. Mullen, Distribution of irreducible polynomials of
small degrees over finite fields, \emph{Math. Comp.} \textbf{67} (1998), 337--341.
\href{https://doi.org/10.1090/S0025-5718-98-00904-1}{doi:10.1090/S0025-5718-98-00904-1}.
\bibitem{HM}
T. Hansen and G. L. Mullen, Primitive polynomials over finite fields,
\emph{Math. Comp.} \textbf{59} (1992), 639--643.
\href{https://doi.org/10.1090/S0025-5718-1992-1134730-7}{doi:10.1090/S0025-5718-1992-1134730-7}.
\bibitem{HIKR}
D. Hart, A. Iosevich, D. Koh and M. Rudnev, Averages over hyperplanes,
sum-product theory in vector spaces over finite fields and the
Erd\H{o}s--Falconer distance conjecture, \emph{Trans. Amer. Math. Soc.}
\textbf{363} (2011), 3255--3275.
\href{https://doi.org/10.1090/S0002-9947-2010-05232-8}{doi:10.1090/S0002-9947-2010-05232-8}.
\bibitem{Hayes}
D. R. Hayes, The expression of a polynomial as a sum of three irreducibles,
\emph{Acta Arith.} \textbf{11} (1966), 461--488.
\href{https://doi.org/10.4064/aa-11-4-461-488}{doi:10.4064/aa-11-4-461-488}.
\bibitem{KapetanakisInvariant}
G. Kapetanakis, Prescribing coefficients of invariant irreducible
polynomials, \emph{J. Number Theory} \textbf{180} (2017), 615--628.
\href{https://doi.org/10.1016/j.jnt.2017.05.003}{doi:10.1016/j.jnt.2017.05.003}.

\bibitem{Merai}
L. M\'erai, Rudin--Shapiro function along irreducible polynomials
over finite fields, \emph{Proc. Amer. Math. Soc.}, to appear.
\href{https://doi.org/10.1090/proc/17975}{doi:10.1090/proc/17975}.
\href{https://arxiv.org/abs/2411.19012v3}{arXiv:2411.19012v3}.

\bibitem{Meyn}
H. Meyn, On the construction of irreducible self-reciprocal polynomials
over finite fields, \emph{Appl. Algebra Engrg. Comm. Comput.}
\textbf{1} (1990), 43--53.
\href{https://doi.org/10.1007/BF01810846}{doi:10.1007/BF01810846}.
\bibitem{Pollack}
P. Pollack, Irreducible polynomials with several prescribed coefficients,
\emph{Finite Fields Appl.} \textbf{22} (2013), 70--78.
\href{https://doi.org/10.1016/j.ffa.2013.03.001}{doi:10.1016/j.ffa.2013.03.001}.
\bibitem{Porritt}
S. Porritt, \emph{Some topics in the analytic number theory of polynomials
over a finite field}, Ph.D. thesis, University College London, 2020.
\url{https://discovery.ucl.ac.uk/id/eprint/10101953/}.
\bibitem{Wan}
D. Wan, Generators and irreducible polynomials over finite fields,
\emph{Math. Comp.} \textbf{66} (1997), 1195--1212.
\href{https://doi.org/10.1090/S0025-5718-97-00835-1}{doi:10.1090/S0025-5718-97-00835-1}.
\bibitem{Weil}
A. Weil, On some exponential sums,
\emph{Proc. Natl. Acad. Sci. USA} \textbf{34} (1948), 204--207.
\href{https://doi.org/10.1073/pnas.34.5.204}{doi:10.1073/pnas.34.5.204}.
\bibitem{YucasMullen}
J. L. Yucas and G. L. Mullen, Self-reciprocal irreducible polynomials
over finite fields, \emph{Des. Codes Cryptogr.} \textbf{33} (2004), 275--281.
\href{https://doi.org/10.1023/B:DESI.0000036251.41345.1f}{doi:10.1023/B:DESI.0000036251.41345.1f}.
\end{thebibliography}
\end{document}